\documentclass[11pt]{article}
\usepackage{mathrsfs}
\usepackage{amsthm}
\usepackage{amssymb}
\usepackage{amsmath}
\usepackage{graphicx}
\usepackage{color}
\usepackage{amsfonts}
\usepackage{float}
\usepackage{cite}
\usepackage[text={140mm,210mm},left=35mm,vmarginratio=1:1]{geometry}
\newtheorem{theorem}{Theorem}[section]
\newtheorem{remark}{Remark}[section]
\newtheorem{lemma}[theorem]{Lemma}

\numberwithin{equation}{section}
\normalsize

\begin{document}
\title{\textbf{Hydrodynamics of $N$-urn susceptible-infected-removed epidemics}}

\author{Xiaofeng Xue \thanks{\textbf{E-mail}: xfxue@bjtu.edu.cn \textbf{Address}: School of Mathematics and Statistics, Beijing Jiaotong University, Beijing 100044, China.}\\ Beijing Jiaotong University}

\date{}
\maketitle

\noindent {\bf Abstract:} In this paper we are concerned with $N$-urn susceptible-infected-removed epidemics, where each urn is in one of three states, namely `susceptible', `infected' and `removed'. We assume that recovery rates of infected urns and infection rates between infected and susceptible urns are
all coordinate-dependent. We show that the hydrodynamic limit of our model is driven by a deterministic $(C[0, 1])^\prime$-valued process with density which is the solution to a nonlinear $C[0, 1]$-valued ordinary differential equation consistent with a mean-field analysis. We further show that the fluctuation of our process is driven by a generalized Ornstein-Uhlenbeck process. A key step in proofs of above main results is to show that states of different urns are approximately independent as $N\rightarrow+\infty$.

\quad


\section{Introduction}\label{section one}
In this paper we are concerned with the $N$-urn susceptible-infected-removed (SIR) model. We first introduce the model. For given integer $N\geq 1$, the $N$-urn SIR model $\{\xi_t\}_{t\geq 0}$ is a continuous-time Markov process with state space
\[
\{-1, 0, 1\}^N=\left\{\xi=\left(\xi(1), \xi(2),\ldots, \xi(N)\right):~\xi(i)\in\{-1, 0, 1\}\text{~for~}1\leq i\leq N\right\}
,
\]
i.e., there are $N$ urns such that at each urn there is a spin taking values in $\{-1, 0, 1\}$. There are two positive functions $\lambda\in C([0, 1]\times[0, 1])$ and $\psi\in C[0, 1]$ such that for all $1\leq i\leq N$:

1) if $\xi_t(i)=1$, then $\xi_t(i)$ jumps to $-1$ at rate $\psi\left(\frac{i}{N}\right)$.

2) If $\xi_t(i)=0$, then $\xi_t(i)$ jumps to $1$ at rate
\[
\frac{1}{N}\sum_{j=1}^N\lambda\left(\frac{i}{N}, \frac{j}{N}\right)1_{\{\xi_t(j)=1\}},
\]
where $1_A$ is the indicator function of the event $A$. We further assume that an urn in state $-1$ is frozen in this state forever. The process $\{\xi_t\}_{t\geq 0}$ can be defined equivalently via its generator. According to the transitions rates given above. The generator $\mathcal{L}$ of $\{\xi_t\}_{t\geq 0}$ is given by
\begin{align}\label{equ 1.1 generator of the SIR}
\mathcal{L}f(\xi)&=\sum_{i=1}^N\psi\left(\frac{i}{N}\right)1_{\{\xi(i)=1\}}\left(f(\xi^{i, -1})-f(\xi)\right) \\
&\text{\quad}+\frac{1}{N}\sum_{i=1}^N\sum_{j=1}^N\lambda\left(\frac{i}{N}, \frac{j}{N}\right)1_{\{\xi(i)=0, \xi(j)=1\}}\left(f(\xi^{i, 1})-f(\xi)\right) \notag
\end{align}
for any $\xi\in \{-1, 0, 1\}^N$ and $f$ from $\{-1, 0, 1\}^N$ to $\mathbb{R}$, where
\[
\xi^{i, k}(j)=
\begin{cases}
\xi(j) & \text{~if~}j\neq i,\\
k  & \text{~if~}j=i
\end{cases}
\]
for all $1\leq i\leq N$ and $k\in\{1, -1\}$. From now on, to distinguish different $N$, we write $\xi_t$ and $\mathcal{L}$ as $\xi_t^N$ and $\mathcal{L}^N$.

Intuitively, $\{\xi_t^N\}_{t\geq 0}$ describes the spread of a SIR epidemic between $N$ individuals. States $-1, 0, 1$ respectively represent `susceptible', `infected' and `removed'. If the $i$th individual is susceptible and the $j$th one is infected, then $j$ infects $i$ at rate $\frac{1}{N}\lambda\left(\frac{i}{N}, \frac{j}{N}\right)$ and $j$ becomes removed at rate $\psi\left(\frac{j}{N}\right)$. A removed individual will never be infected again.

When $\psi\equiv 1$ and $\lambda\equiv \lambda_0\in \mathbb{R}$, our model reduces to the classic spatial-homogeneous SIR model (see \cite{Anderson1991}) such that $\left(\mathcal{S}_t^N, \mathcal{I}_t^N\right)$ jumps to $\left(\mathcal{S}_t^N, \mathcal{I}_t^N-1\right)$ at rate $\mathcal{I}_t^N$ or jumps to $\left(\mathcal{S}_t^N-1, \mathcal{I}_t^N+1\right)$ at rate $\frac{\lambda_0}{N}\mathcal{S}_t^N\mathcal{I}_t^N$, where
\[
\mathcal{S}_t^N=\sum_{i=1}^N1_{\{\xi_t^N(i)=0\}} \text{~and~}\mathcal{I}_t^N=\sum_{i=1}^N1_{\{\xi_t^N(i)=1\}}
\]
are total numbers of susceptible individuals and infected individuals respectively. For the classic SIR model, the law of large numbers (LLN) and the central limit theorem (CLT) of $\{\left(\mathcal{S}_t^N, \mathcal{I}_t^N\right)\}_{t\geq 0}$ can be derived directly as an application of the theory of density-dependent Markov chains given in \cite{Kurtz1978} and Chapter 11 of \cite{Ethier1986}. For a detailed recall, see Subsection \ref{subsection 2.2}.

In this paper, we will investigate the hydrodynamic limit and fluctuation of our process $\{\xi_t^N\}_{t\geq 0}$, i.e., LLN and CLT of
\[
\frac{1}{N}\sum_{i=1}^N1_{\{\xi_t^N(i)=0\}}f\left(\frac{i}{N}\right) \text{~and~} \frac{1}{N}\sum_{i=1}^N1_{\{\xi_t^N(i)=1\}}f\left(\frac{i}{N}\right)
\]
for any $f\in C[0, 1]$ as $N\rightarrow+\infty$. This investigation is motivated by the aforesaid LLN and CLT given in \cite{Kurtz1978} of $\mathcal{S}_t^N$ and $\mathcal{I}_t^N$ in the spatial-homogenous case, where $f\equiv 1$ and recovery and infection rates do not rely on coordinates of individuals. So it is natural to discuss general cases where recovery and infection rates are coordinate-dependent and $f$ can be anyone in $C[0, 1]$. Since we discuss all $f\in C[0, 1]$, we consider empirical density fields
\[
\frac{1}{N}\sum_{i=1}^N1_{\{\xi_t(i)=0\}}\delta_{\frac{i}{N}}(du) \text{~and~} \frac{1}{N}\sum_{i=1}^N1_{\{\xi_t(i)=1\}}\delta_{\frac{i}{N}}(du)
\]
as linear operators from $C[0, 1]$ to $\mathbb{R}$ and will directly give LLN and CLT of these $(C[0, 1])^\prime$-valued random elements. For mathematical details, see Section \ref{section two}.

We are also inspired by previous investigations about hydrodynamics of different types of interacting particle systems on $N$ urns such as Ehrenfest models, voter models, symmetric exclusion processes and binary contact path processes (see \cite{Xue2022} and \cite{Xue2022b}). For those models, hydrodynamic limits are driven by $(C[0, 1])^\prime$-valued linear ordinary differential equations (ODE) and fluctuations are driven by $(C[0, 1])^\prime$-valued Ornstein-Uhlenbeck (O-U) processes. It is natural to ask whether analogue results hold for the SIR model. We will show that the hydrodynamic limit of our SIR model is driven by a deterministic $(C[0, 1])^\prime$-valued process with density which is the solution to a nonlinear $C[0, 1]$-valued ODE and the fluctuation of our SIR model is driven by a $(C[0, 1])^\prime$-valued O-U process. For mathematical details, see Section \ref{section two}.

We are also inspired by the investigation of mean-field limits of long range contact processes (SIS epidemics) on lattices given in \cite{Swindle}. Main results of this paper can be considered as analogues of theorems in \cite{Swindle} since hydrodynamic limits of both models are consist with mean-field analyses and key steps in proofs of respective theorems are both verifications of approximate independence between states of different vertices. However, above verifications of respective models are different in techniques. For the contact process in \cite{Swindle}, the above approximate independence is derived according to the self-duality of the model (see Section 3.4 of \cite{Lig1985}). Since the SIR model does not have a duality property, we adopt a different approach where we investigate infection paths ending at different vertices before a given moment and show that they are not intersected with high probability as $N\rightarrow+\infty$. For mathematical details, see Section \ref{section three}. 

\section{Main results}\label{section two}
In this section, we give our main results and apply them in the spatial-homogeneous case to revisit the LLN and CLT given in \cite{Kurtz1978}.

\subsection{Main results}\label{subsection 2.1}
In this subsection, we give our main result. For later use, we first introduce some notations and definitions. For any $t\geq 0$ and $1\leq i\leq N$, we define
\[
S_t^N(i)=1_{\{\xi_t^N(i)=0\}} \text{~and~} I_t^N(i)=1_{\{\xi_t^N(i)=1\}}.
\]
For any $a\in[0, 1]$, let $\delta_a(du)$ be the Dirac measure concentrated on $a$, then we define
\[
\mu_t^N(du)=\frac{1}{N}\sum_{i=1}^NI_t^N(i)\delta_{\frac{i}{N}}(du)
\]
and
\[
\theta_t^N(du)=\frac{1}{N}\sum_{i=1}^NS_t^N(i)\delta_{\frac{i}{N}}(du)
\]
for all $t\geq 0$ and $N\geq 1$. Consequently, $\mu_t^N$ and $\theta_t^N$ can be considered as random elements in $(C[0, 1])^\prime$ such that
\[
\mu_t^N(f)=\frac{1}{N}\sum_{i=1}^NI_t^N(i)f\left(\frac{i}{N}\right) \text{~and~}\theta_t^N(f)=\frac{1}{N}\sum_{i=1}^NS_t^N(i)f\left(\frac{i}{N}\right)
\]
for any $f\in C[0, 1]$. Throughout this paper, we adopt the following initial assumption.

\textbf{Assumption} (A): $\{\xi_0^N(i)\}_{i=1}^N$ are independent and there exists $\phi\in C[0, 1]$ such that
\[
P\left(\xi_0^N(i)=1\right)=\phi\left(\frac{i}{N}\right)=1-P\left(\xi_0^N(i)=0\right)
\]
for all $1\leq i\leq N$.

The following theorem is our first main result, which gives LLNs of $\mu_t^N$ and $\theta_t^N$ as $N\rightarrow+\infty$, i.e., the hydrodynamic limit of $\xi_t^N$.

\begin{theorem}\label{theorem 2.1 hydrodynamic limit}
Under Assumption (A),
\[
\lim_{N\rightarrow+\infty}\mu_t^N(f)=\int_0^1 \rho_1(t,u)f(u)du \text{~and~}\lim_{N\rightarrow+\infty}\theta_t^N(f)=\int_0^1\rho_0(t,u)f(u)du
\]
in probability for any $f\in C[0, 1]$, where $\left(\rho_1(t, \cdot), \rho_0(t, \cdot)\right)$ is the unique solution to the $(C[0,1])^2$-valued ODE
\begin{equation}\label{equ 2.1 nonliear ODE}
\begin{cases}
&\frac{d}{dt}\rho_1(t,u)=-\psi(u)\rho_1(t,u)+\rho_0(t,u)\int_0^1\lambda(u,v)\rho_1(t,v)dv, \text{~}0\leq u\leq 1,\\
&\frac{d}{dt}\rho_0(t,u)=-\rho_0(t,u)\int_0^1\lambda(u, v)\rho_1(t,v)dv, \text{~}0\leq u\leq 1,\\
&\rho_1(0, u)=\phi(u), \text{~}0\leq u\leq 1,\\
&\rho_0(0, u)=1-\phi(u), \text{~}0\leq u\leq 1.
\end{cases}
\end{equation}
\end{theorem}

\begin{remark}\label{remark 2.1}
Note that Equation \eqref{equ 2.1 nonliear ODE} satisfies local Lipschitz's condition under the norm $\|(\cdot, \cdot)\|_\infty$ on $(C[0, 1])^2$ given by
\[
\|(f,g)\|_\infty=\sup_{0\leq u\leq 1}|f(u)|+\sup_{0\leq u\leq 1}|g(u)|
\]
for all $(f,g)\in (C[0, 1])^2$, hence the existence and uniqueness of the solution to Equation \eqref{equ 2.1 nonliear ODE} for $t$ in a local area follows from Theorem 1.2 in Chapter 19 of \cite{Lang}. Furthermore, it is easy to check that the solution $\left(\rho_1(t, \cdot), \rho_0(t, \cdot)\right)$ satisfies that
$\|\left(\rho_1(t, \cdot), \rho_0(t, \cdot)\right)\|_\infty\leq 1$ and hence the solution to Equation \eqref{equ 2.1 nonliear ODE} can be defined for $t\in [0, +\infty)$.
\end{remark}

\begin{remark}\label{remark 2.2}
A mean field analysis gives an intuitive explanation of Theorem \ref{theorem 2.1 hydrodynamic limit}. Our SIR model has the Kolmogorov-Chapman (KC) equation
\begin{align*}
\frac{d}{dt}P(I_t(i)=1)=&-\psi\left(\frac{i}{N}\right)P(I_t(i)=1)\\
&+\frac{1}{N}\sum_{j=1}^N\lambda\left(\frac{i}{N}, \frac{j}{N}\right)P\left(S_t(i)=1, I_t(j)=1\right).
\end{align*}
For large $N$, since the infection rate between two different urns is $O(\frac{1}{N})$, it is natural to think that states of different urns are approximately independent and hence the term $P\left(S_t(i)=1, I_t(j)=1\right)$ in the KC equation can be replaced by
\[
P\left(S_t(i)=1\right)P\left(I_t(j)=1\right)
\]
with small error, which gives a discretely approximate version of Equation \eqref{equ 2.1 nonliear ODE}. In fact, the idea of our proof of Theorem \ref{theorem 2.1 hydrodynamic limit} is based on the above mean field analysis. We will rigorously prove that ${\rm Cov}(S_t^N(i), I_t^N(j))$ converges to $0$ as $N\rightarrow+\infty$ when $i\neq j$. For mathematical details, see Section \ref{section three}.
\end{remark}

Now we discuss the fluctuation of our process $\{\xi_t^N\}_{t\geq 0}$, i.e., the CLT from the LLN given in Theorem \ref{theorem 2.1 hydrodynamic limit}. To give our main result, we introduce some nations and definitions. For $N\geq 1$ and $t\geq 0$, we define
\[
\eta_t^N(du)=\frac{1}{\sqrt{N}}\sum_{i=1}^N\left(I_t^N(i)-\mathbb{E}I_t^N(i)\right)\delta_{\frac{i}{N}}(du)
\]
and
\[
\beta_t^N(du)=\frac{1}{\sqrt{N}}\sum_{i=1}^N\left(S_t^N(i)-\mathbb{E}S_t^N(i)\right)\delta_{\frac{i}{N}}(du),
\]
where $\mathbb{E}$ is the expectation operator. Let $(C[0, 1])^\prime$ be equipped with the weak topology, i.e., $\nu_n$ converges to $\nu$ in $(C[0, 1])^\prime$ if and only if
\[
\lim_{n\rightarrow+\infty}\nu_n(f)=\nu(f)
\]
for all $f\in C[0, 1]$, then $\{\eta_t^N\}_{0\leq t\leq T}$ and $\{\beta_t^N\}_{0\leq t\leq T}$ can be considered as random elements in $\mathcal{D}\left([0, T], (C[0, 1])^\prime\right)$
for given $T>0$, where $\mathcal{D}\left([0, T], (C[0, 1])^\prime\right)$ is the set of c\`{a}dl\`{a}g functions from $[0, T]$ to $(C[0, 1])^\prime$ equipped with the Skorohod topology.

We denote by $\{\mathcal{W}_t\}_{t\geq 0}$ the $(C[0, 1])^\prime$-valued standard Brownian motion such that $\{\mathcal{W}_t(f)\}_{t\geq 0}$ is a real-valued Brownian motion with
\[
{\rm Cov}\left(\mathcal{W}_t(f), \mathcal{W}_t(f)\right)=t\int_0^1f^2(u)du
\]
for any $f\in C[0, 1]$ and $t\geq 0$. For any $t\geq 0$, we denote by $\mathcal{A}_{1,t}$ the linear operator from $C[0, 1]$ to $C[0, 1]$ such that
\[
\mathcal{A}_{1,t}f(u)=f(u)\int_0^1\lambda(u,v)\rho_1(t,v)dv
\]
and by $\mathcal{A}_{0, t}$ the linear operator from $C[0, 1]$ to $C[0, 1]$ such that
\[
\mathcal{A}_{0, t}f(u)=\int_0^1\lambda(v,u)\rho_0(t,v)f(v)dv
\]
for any $f\in C[0,1]$ and $u\in [0, 1]$, where $\left(\rho_1(t,\cdot), \rho_0(t,\cdot)\right)$ is the solution to Equation \eqref{equ 2.1 nonliear ODE}. We denote by $\mathcal{A}_\psi$ the linear operator from $C[0, 1]$ to $C[0, 1]$ such that
\[
\mathcal{A}_\psi f(u)=\psi(u)f(u)
\]
for any $f\in C[0, 1], u\in [0, 1]$. For any $t\geq 0$, we denote by $b_t$ the linear operator from $C[0, 1]$ to $C[0, 1]$ such that
\[
b_tf(u)=\sqrt{\psi(u)\rho_1(t,u)}f(u)
\]
and by $\alpha_t$ the linear operator from $C[0, 1]$ to $C[0, 1]$ such that
\[
\alpha_tf(u)=\sqrt{\int_0^1\lambda(u,v)\rho_0(t, u)\rho_1(t,v)dv}f(u)
\]
for any $f\in C[0, 1]$ and $u\in [0, 1]$. For any linear operator $\mathcal{P}$ from $C[0, 1]$ to $C[0, 1]$, we denote by $\mathcal{P}^{*}$ the adjoint operator of $\mathcal{P}$, i.e., $\mathcal{P}^{*}$ is the linear operator from $(C[0, 1])^\prime$ to $(C[0, 1])^\prime$ such that
\[
\mathcal{P}^{*}\nu(f)=\nu(\mathcal{P}f)
\]
for any $f\in C[0, 1]$ and $\nu\in (C[0, 1])^\prime$.

Let $\{\mathcal{W}_t^1\}_{t\geq 0}$ and $\{\mathcal{W}_t^2\}_{t\geq 0}$ be two independent copies of $\{\mathcal{W}_t\}_{t\geq 0}$, then we introduce the definition of the solution to the following $((C[0, 1])^\prime)^2$-valued stochastic differential equation (SDE):
\begin{equation}\label{equ 2.2 generalized O-U process}
\begin{cases}
&d\eta_t=\mathcal{A}_{1,t}^{*}\beta_tdt+\mathcal{A}_{0, t}^{*}\eta_tdt-\mathcal{A}_\psi^{*}\eta_tdt+b_t^*d\mathcal{W}_t^1+\alpha_t^*d\mathcal{W}_t^2,\\
&d\beta_t=-\mathcal{A}_{1,t}^{*}\beta_tdt-\mathcal{A}_{0, t}^{*}\eta_tdt-\alpha_t^*d\mathcal{W}_t^2,\\
&(\eta_0, \beta_0)\text{~is independent of~}\{\mathcal{W}_t^1\}_{t\geq 0}\text{~and~}\{\mathcal{W}_t^2\}_{t\geq 0}.
\end{cases}
\end{equation}
According to an analysis similar with that leading to Theorem 1.4 of \cite{Holley1978}, with given initial distribution of $(\eta_0, \beta_0)$, there exists an unique random element $\{(\eta_t, \beta_t)\}_{t\geq 0}$ such that:

1) $\{(\eta_t(f), \beta_t(f))\}_{t\geq 0}$ is a $\mathbb{R}^2$-valued continuous process for any $f\in C[0, 1]$.

2) For any $G\in C_c^\infty(\mathbb{R}^2)$ and $f, g\in C[0, 1]$,
\begin{align*}
\Bigg\{&G(\eta_t(f), \beta_t(g))-G(\eta_0(f), \beta_0(g))\\
&\text{~}-\int_0^t\frac{\partial}{\partial u}G(\eta_s(f), \beta_s(g))\left(\beta_s(\mathcal{A}_{1,s}f)+\eta_s(\mathcal{A}_{0,s}f-\mathcal{A}_\psi f)\right)ds\\
&\text{~}+\int_0^t\frac{\partial}{\partial v}G(\eta_s(f), \beta_s(g))\left(\beta_s(\mathcal{A}_{1,s}g)+\eta_s(\mathcal{A}_{0, s}g)\right)ds\\
&\text{~}-\frac{1}{2}\int_0^t\frac{\partial^2}{\partial u^2}G(\eta_s(f), \beta_s(g))\\
&\text{\quad\quad\quad}\times\left(\int_0^1\rho_1(s,u)\psi(u)f^2(u)du+\int_0^1\int_0^1\lambda(u,v)\rho_0(s,u)\rho_1(s,v)f^2(u)dudv\right)ds\\
&\text{~}-\frac{1}{2}\int_0^t\int_0^1\int_0^1\frac{\partial^2}{\partial v^2}G(\eta_s(f), \beta_s(g))\lambda(u,v)\rho_0(s,u)\rho_1(s,v)g^2(u)dudvds\\
&\text{~}+\int_0^t\int_0^1\int_0^1\frac{\partial^2}{\partial uv}G(\eta_s(f), \beta_s(g))\lambda(u,v)\rho_0(s,u)\rho_1(s,v)f(u)g(u)dudvds\Bigg\}_{t\geq 0}
\end{align*}
is a martingale, where $u$ is the first coordinate and $v$ is the second coordinate. According to It\^{o}'s formula, it is nature to define the above $\{(\eta_t, \beta_t)\}_{t\geq 0}$ as the solution to Equation \eqref{equ 2.2 generalized O-U process}, which is called a generalized O-U process as introduced in \cite{Holley1978}.

\begin{remark}\label{remark 2.3}
As introduced in \cite{Holley1978}, the formal definition of the solution to Equation \eqref{equ 2.2 generalized O-U process} is the unique solution to the corresponding martingale problem. However, the analysis leading to Theorem 1.4 of \cite{Holley1978} implies that we can define the solution to Equation \eqref{equ 2.2 generalized O-U process} equivalently through solving this equation directly. In detail, Equation \eqref{equ 2.2 generalized O-U process} can be written as
\[
d\begin{pmatrix} \eta_t\\ \beta_t\end{pmatrix}=\Sigma_{1,t}^*\begin{pmatrix} \eta_t\\ \beta_t\end{pmatrix}dt+\Sigma_{2,t}^*d\begin{pmatrix} \mathcal{W}_t^1\\ \mathcal{W}_t^2\end{pmatrix},
\]
where
\[
\Sigma_{1,t}^*=\begin{pmatrix} \mathcal{A}_{1,t}^*&\mathcal{A}_{0,t}^*-\mathcal{A}_\psi^*\\-\mathcal{A}_{1,t}^*&-\mathcal{A}_{0,t}^*\end{pmatrix}
\text{~and~}\Sigma_{2,t}^*=\begin{pmatrix}b_t^*&\alpha_t^*\\0&-\alpha_t^*\end{pmatrix}.
\]
Let $\{\Upsilon^*(t,s)\}_{0\leq s\leq t<+\infty}$ be the solution to the matrix-valued ODE
\begin{equation}\label{equ 2.3 auxiliary ODE}
\begin{cases}
&\frac{\partial}{\partial t}\Upsilon^*(t, s)=\Sigma_{1,t}^*\Upsilon^*(t,s),\\
&\Upsilon^*(s,s)=\begin{pmatrix}I_d^*&0 \\0& I_d^*\end{pmatrix}, s\geq 0,
\end{cases}
\end{equation}
where $I_d^*$ is the identity function on $(C[0, 1])^\prime$. Note that the existence and uniqueness of $\Upsilon^*$ follows from Theorem 1.2 in Chapter 19 of \cite{Lang}. Then, according to It\^{o}'s formula,
\begin{equation}\label{equ 2.4 strong solution}
\begin{pmatrix} \eta_t\\ \beta_t\end{pmatrix}=\Upsilon^*(t,0)\begin{pmatrix} \eta_0\\ \beta_0\end{pmatrix}+\int_0^t\Upsilon^*(t,s)\Sigma_{2,s}^*d\begin{pmatrix} \mathcal{W}_s^1\\ \mathcal{W}_s^2\end{pmatrix}.
\end{equation}
\end{remark}

Now we can state our second main result, which gives the fluctuation of our SIR model.
\begin{theorem}\label{theorem 2.2 fluctuation}
Under Assumption (A), for given $T>0$, $\{(\eta_t^N, \beta_t^N)\}_{0\leq t\leq T}$ converges weakly to $\{(\eta_t, \beta_t)\}_{0\leq t\leq T}$ as $N\rightarrow+\infty$, where $\{(\eta_t, \beta_t)\}_{t\geq 0}$ is the solution to Equation \eqref{equ 2.2 generalized O-U process} with the following initial distribution:

1) $\beta_0=-\eta_0$.

2) For any $f\in C[0, 1]$, $\eta_0(f)$ follows $\mathbb{N}(0, \int_0^1f^2(u)\phi(u)(1-\phi(u))du)$, where $\mathbb{N}(a, \sigma^2)$ is the normal distribution with mean $a$ and variance $\sigma^2$.
\end{theorem}

Theorem \ref{theorem 2.2 fluctuation} shows that the fluctuation of $\{\xi_t^N\}_{t\geq 0}$ is driven by a generalized O-U process. According to Remark \ref{remark 2.3}, for any $f,g\in C[0, 1]$ and $t, s\geq 0$, $(\eta_t^N(f), \beta_s^N(g))$ converges weakly to the normal distribution of $(\eta_t(f), \beta_s(g))$ given by \eqref{equ 2.4 strong solution}.

Remaining sections of this paper are organized as follows. We prove a replacement lemma in Section \ref{section three} as a preliminary, which shows that states of different urns are approximately independent as $N\rightarrow+\infty$. We prove Theorem \ref{theorem 2.1 hydrodynamic limit} in Section \ref{section four}. With our replacement lemma, Theorem \ref{theorem 2.1 hydrodynamic limit} follows from a classic mean-variance analysis. We prove Theorem \ref{theorem 2.2 fluctuation} in Section \ref{section five}. According to the Dynkin's martingale formula and our replacement lemma, we can show that each weak limit of subsequence of $\{(\eta_t^N, \beta_t^N):~0\leq t\leq T\}_{N\geq 1}$ is the solution to the martingale problem with respect to Equation \eqref{equ 2.2 generalized O-U process} and hence Theorem \ref{theorem 2.2 fluctuation} holds after we check that $\{(\eta_t^N, \beta_t^N):~0\leq t\leq T\}_{N\geq 1}$ are tight.

\subsection{A revisit of the spatial-homogeneous case} \label{subsection 2.2}

As an application of our main result, in this section we revisit LLN and CLT of the spatial-homogeneous SIR model given in \cite{Kurtz1978}, i.e., throughout this subsection we assume that $\psi\equiv 1$, $\lambda\equiv \lambda_0\in \mathbb{R}$ and $\phi\equiv \phi_0\in [0, 1]$. As we have introduced in Section \ref{section one}, $\{(\mathcal{I}_t^N, \mathcal{S}_t^N)^{\mathsf{T}}\}_{t\geq 0}$ is a $\mathbb{R}^2$-valued continuous-time Markov process such that
\[
(\mathcal{I}_t^N, \mathcal{S}_t^N)^{\mathsf{T}} \text{~jumps to~} (\mathcal{I}_t^N, \mathcal{S}_t^N)^{\mathsf{T}}+l_k \text{~at rate~}NF_k\left(\frac{\mathcal{I}_t^N}{N}, \frac{\mathcal{S}_t^N}{N}\right)
\]
for $k=1, 2$, where $\mathsf{T}$ is the transposition operator, $l_1=(-1, 0)^{\mathsf{T}}, l_2=(1, -1)^{\mathsf{T}}$, $F_1(u,v)=u$ and $F_2(u,v)=\lambda_0uv$. As a result, $\{(\mathcal{I}_t^N, \mathcal{S}_t^N)^{\mathsf{T}}\}_{t\geq 0}$ is an example of the density-dependent Markov chain (DDMC) introduced in \cite{Kurtz1978}. Then, as direct applications of main results given in \cite{Kurtz1978}, we have the following two theorems.

\begin{theorem}\label{theorem 2.3 LLN homogeneous case}(Kurtz, 1978)
As $N\rightarrow+\infty$, $(\mathcal{I}_t^N/N, \mathcal{S}_t^N/N)^{\mathsf{T}}$ converges in probability to $x_t$, where
\begin{equation}\label{equ 2.5 LLN homogeneous}
\begin{cases}
&\frac{d}{dt}x_t=\sum_{k=1}^2l_kF_k(x_t),\\
&x_0=(\phi_0, 1-\phi_0)^{\mathsf{T}}.
\end{cases}
\end{equation}
\end{theorem}

\begin{theorem}\label{theorem 2.4 CLT homogeneous case}(Kurtz, 1978)
As $N\rightarrow+\infty$, $\left(\frac{\mathcal{I}_t^N-\mathbb{E}\mathcal{I}_t^N}{\sqrt{N}}, \frac{\mathcal{S}_t^N-\mathbb{E}\mathcal{S}_t^N}{\sqrt{N}}\right)^{\mathsf{T}}$ converges in distribution to $V_t$, where
\begin{equation}\label{equ 2.6 CLT homogeneous}
dV_t=\sum_{k=1}^2l_k\nabla^\mathsf{T}F_k(x_t)V_tdt+\left(\sum_{k=1}^2l_kF_k(x_t)l_k^\mathsf{T}\right)^{\frac{1}{2}}d\mathcal{B}_t,
\end{equation}
where $\mathcal{B}_t$ is a $2$-dimensional standard Brownian motion and $\nabla^\mathsf{T}=(\frac{\partial}{\partial u}, \frac{\partial}{\partial v})$.
\end{theorem}

In \cite{Kurtz1978}, LLN and CLT of DDMC are derived through writing the process as the solution to a SDE driven by Poisson processes. Here we prove Theorems \ref{theorem 2.3 LLN homogeneous case} and \ref{theorem 2.4 CLT homogeneous case} via our main results Theorems \ref{theorem 2.1 hydrodynamic limit} and \ref{theorem 2.2 fluctuation}.

\proof[An alternative proof of Theorem \ref{theorem 2.3 LLN homogeneous case}]

Since $\frac{1}{N}\mathcal{I}_t^N=\mu_t^N(1)$ and $\frac{1}{N}\mathcal{S}_t^N=\theta_t^N(1)$, by Theorem \ref{theorem 2.1 hydrodynamic limit}, $\frac{1}{N}\mathcal{I}_t^N, \frac{1}{N}\mathcal{S}_t^N$ converges in probability to
\[
i_t=\int_0^1\rho_1(t,u)du \text{~and~}s_t=\int_0^1\rho_0(t,u)du
\]
respectively, where
\[
\frac{d}{dt}\rho_1(t,u)=-\rho_1(t,u)+\lambda_0\rho_0(t,u)\int_0^1\rho_1(t,v)dv,
\]
and
\[
\frac{d}{dt}\rho_0(t,u)=-\lambda_0\rho_0(t,u)\int_0^1\rho_1(t,v)dv.
\]
As a result,
\[
\begin{cases}
&\frac{d}{dt}i_t=-i_t+\lambda_0i_ts_t,\\
&\frac{d}{dt}s_t=-\lambda_0i_ts_t,\\
&(i_0, s_0)=(\phi_0, 1-\phi_0),
\end{cases}
\]
which is Equation \eqref{equ 2.5 LLN homogeneous}.

\qed

\proof[An alternative proof of Theorem \ref{theorem 2.4 CLT homogeneous case}]

Since
\[
\frac{\mathcal{I}_t^N-\mathbb{E}\mathcal{I}_t^N}{\sqrt{N}}=\eta_t^N(1) \text{~and~} \frac{\mathcal{S}_t^N-\mathbb{E}\mathcal{S}_t^N}{\sqrt{N}}=\beta_t^N(1),
\]
by Theorem \ref{theorem 2.2 fluctuation}, $\left(\frac{\mathcal{I}_t^N-\mathbb{E}\mathcal{I}_t^N}{\sqrt{N}}, \frac{\mathcal{S}_t^N-\mathbb{E}\mathcal{S}_t^N}{\sqrt{N}}\right)^{\mathsf{T}}$ converges in distribution to $\left(\eta_t(1), \beta_t(1)\right)^\mathsf{T}$,
where
\[
\begin{cases}
&d\eta_t(1)=\beta_t(\mathcal{A}_{1,t}1)dt+\eta_t(\mathcal{A}_{0, t}1)dt-\eta_t(\mathcal{A}_{\psi}1)dt+\|b_t1\|_2d\hat{B}_t^1+\|\alpha_t1\|_2d\hat{B}_t^2,\\
&d\beta_t(1)=-\beta_t(\mathcal{A}_{1,t}1)dt-\eta_t(\mathcal{A}_{0,t}1)dt-\|\alpha_t1\|_2d\hat{B}_t^2,
\end{cases}
\]
where $\|f\|_2=\sqrt{\int_0^1f^2(u)du}$ for any $f\in C[0, 1]$ and $\{(\hat{B}_t^1, \hat{B}_t^2)\}_{t\geq 0}$ is a $\mathbb{R}^2$-valued standard Brownian motion. By direct calculation,
\[
\beta_t(\mathcal{A}_{1,t}1)=\lambda_0i_t\beta_t(1), \text{~}\eta_t(\mathcal{A}_{0, t}1)=\lambda_0s_t\eta_t(1) \text{~and~}\mathcal{A}_\psi1=1.
\]
Therefore, let $V_t=\left(\eta_t(1), \beta_t(1)\right)^\mathsf{T}$, then
\[
dV_t=\begin{pmatrix} \lambda_0s_t-1& \lambda_0i_t\\-\lambda_0s_t&-\lambda_0i_t\end{pmatrix}V_tdt
+\begin{pmatrix} \|b_t1\|_2& \|\alpha_t1\|_2\\0&-\|\alpha_t1\|_2\end{pmatrix}d\begin{pmatrix} \hat{B}_t^1\\ \hat{B}_t^2\end{pmatrix}.
\]
By direct calculation, it is easy to check that
\[
\sum_{k=1}^2l_k\nabla^\mathsf{T}F_k(x_t)=\begin{pmatrix} \lambda_0s_t-1& \lambda_0i_t\\-\lambda_0s_t&-\lambda_0i_t\end{pmatrix}
\]
and
\[
\begin{pmatrix} \|b_t1\|_2& \|\alpha_t1\|_2\\0&-\|\alpha_t1\|_2\end{pmatrix}\begin{pmatrix} \|b_t1\|_2& \|\alpha_t1\|_2\\0&-\|\alpha_t1\|_2\end{pmatrix}^\mathsf{T}=\sum_{k=1}^2l_kF_k(x_t)l_k^\mathsf{T}
\]
and hence $\{V_t\}_{t\geq 0}$ is the solution to Equation \eqref{equ 2.6 CLT homogeneous}.

\qed

The large deviation principle of $\left(\mathcal{I}_t^N, \mathcal{S}_t^N\right)$ for the classis SIR model is given in \cite{Pardoux2017}, so it is natural to ask whether there is a large deviation principle for our process $\left(\mu_t^N, \theta_t^N\right)$ through which and the contraction principle we can revisit the theorem in \cite{Pardoux2017}. We will work on this question as a further investigation.

\section{Replacement lemma}\label{section three}

In this section we prove a replacement lemma which shows that states of different urns are approximately independent as $N\rightarrow+\infty$. For simplicity, we define
\[
H_{t, 0}^N(i)=S_t^N(i) \text{~and~} H_{t, 1}^N(i)=I_t^N(i)
\]
for $1\leq i\leq N$, then our replacement lemma is as follows.

\begin{lemma}\label{lemma 3.1 replacement}
There exists $C_1, C_2, C_3<+\infty$ independent of $N$ such that for all $0\leq t\leq T$, $N\geq 4$, $(r_1, r_2, r_3, r_4)\in \{0, 1\}^4$ and any four different $i,j,k,l\in \{1,2,\ldots,N\}$,
\begin{align}\label{equ 3.1 replacement lemma partOne}
\left|\mathbb{E}(H_{t, r_1}^N(i)H_{t, r_2}^N(j))-\mathbb{E}H_{t, r_1}^N(i)\mathbb{E}H_{t, r_2}^N(j)\right|\leq \frac{C_1}{N},
\end{align}
\begin{align}\label{equ 3.1 replacement lemma partTwo}
\left|\mathbb{E}(H_{t, r_1}^N(i)H_{t, r_2}^N(j)H_{t, r_3}^N(k))-\mathbb{E}H_{t, r_1}^N(i)\mathbb{E}H_{t, r_2}^N(j)\mathbb{E}H_{t, r_3}^N(k)\right|\leq \frac{C_2}{N}
\end{align}
and
\begin{align}\label{equ 3.1 replacement lemma partThree}
\left|\mathbb{E}(H_{t, r_1}^N(i)H_{t, r_2}^N(j)H_{t, r_3}^N(k)H_{t, r_4}^N(l))-\mathbb{E}H_{t, r_1}^N(i)\mathbb{E}H_{t, r_2}^N(j)\mathbb{E}H_{t, r_3}^N(k)\mathbb{E}H_{t, r_4}^N(l)\right|\leq \frac{C_3}{N}.
\end{align}
\end{lemma}

The idea of the proof of Lemma \ref{lemma 3.1 replacement} is as follows.  We introduce some independent exponential times describing recovery times of infected vertices and infection times between infected and susceptible vertices. Then, $\{\xi_t\}_{t>0}$ can be defined as being measurable with respect to the sigma-algebra generated by $\xi_0$ and above exponential times. By further introducing some independent copies of above exponential times, we can construct $\hat{\xi}_t(j), \hat{\xi}_t(k), \hat{\xi}_t(l)$ such that:

1) $\left(\xi_t(i), \hat{\xi}_t(j), \hat{\xi}_t(k), \hat{\xi}_t(l)\right)$ are independent.

2) $\hat{\xi}_t(m)$ and $\xi_t(m)$ have the same distribution for $m=j, k, l$.

3) Conditioned on an event with probability $1-O(1/N)$,
\[
\left(\hat{\xi}_t(j), \hat{\xi}_t(k), \hat{\xi}_t(l)\right)=\left(\xi_t(j), \xi_t(k), \xi_t(l)\right).
\]

In detail, for each $1\leq i\leq N$, let $K_i^N$ be an exponential time with rate $\psi\left(\frac{i}{N}\right)$. For any $1\leq i\neq j\leq N$, let $U_{(i, j)}^N$ be an exponential time with rate $\frac{1}{N}\lambda\left(\frac{i}{N}, \frac{j}{N}\right)$. We further assume that all these exponential times are independent. Meanings of these exponential times are as follows. For any $1\leq i\neq j\leq N$, if $j$ is infected at some moment, then $j$ waits for $K_j^N$ to become removed and waits for $U_{(i, j)}^N$ to infect $i$. This infection really occurs when $U_{(i, j)}^N<K_j^N$ and $i$ is susceptible at the moment just before being infected by $j$. According to the above explanation, let $\xi_0^N$ be distributed as in Assumption (A), then $\xi_t^N(i)$ is determined by $\xi_0^N$ and above exponential times. Under Assumption (A), $\xi_s^N(i)=1$ for some $0\leq s\leq t$ when and only when there exist integer $n\geq 0$ and self-avoiding path $m_0, m_1, m_2,..., m_n=i$ on $\{1,2,\ldots,N\}$ such that

1) $\xi_0^N(m_0)=1$.

2) $U_{(m_{l+1}, m_l)}^N<K_{m_l}^N$ for $0\leq l\leq n-1$.

3) $\sum_{l=0}^{n-1}U^N_{(m_{l+1}, m_l)}\leq t$.

For $i$ has no self-avoiding path with above properties, $\xi_t^N(i)=0$. For $i$ has aforesaid self-avoiding paths, we choose the path making $\sum_{l=0}^{n-1}U^N_{(m_{l+1}, m_l)}$ minimum. Let $c$ be the aforesaid minimum, i.e., $c$ is the moment at which $i$ is infected. Then $\xi_t^N(i)=1$ when $c+K_{i}^N>t$ and $\xi_t^N(i)=-1$ when $c+K_i^N\leq t$. The above description of the SIR model can be considered as an application of the graphical method introduced in \cite{Har1978}.

To define copies of $\xi_t^N(j), \xi_t^N(k), \xi_t^N(l)$, we introduce following notations. For any $1\leq m\neq q\leq N$, let $\{\xi_{0, r}^N(m)\}_{r=2,3,4}$ be three independent copies of $\xi_0^N(m)$, $\{K_{m,r}^N\}_{r=2,3,4}$ be three independent copies of $K_m^N$ and $\{U_{(m,q),r}^N\}_{r=2,3,4}$ be three independent copies of $U_{(m,q)}^N$. We assume all these random variables are independent. For $t>0, 1\leq m\leq N$ and $0\leq r\leq N-1$, we define random set $\Xi^{N, t}_{m,r}\subseteq \{1,2,\ldots,N\}$ by induction as follows.

1) $\Xi^{N, t}_{m, 0}=\{m\}$.

2) After $\Xi^{N,t}_{m, r}$ are defined for $0\leq r\leq q<N-1$, $\Xi^{N, t}_{m, q+1}$ is defined as the set of $\varphi\in \{1,2,\ldots, N\}\setminus\left(\bigcup_{r=0}^q\Xi^{N, t}_{m, r}\right)$ such that there exists self-avoiding path $\varphi=m_{q+1}, m_q, \ldots, m_1, m_0=m$ such that $m_r\in \bigcup_{l=1}^r\Xi^{N, t}_{m,l}$ for $1\leq r\leq q$ and
\[
\sum_{r=0}^{q}U^N_{(m_r, m_{r+1})}\leq t.
\]
Note that if $\Xi^{N, t}_{m, r}=\emptyset$ for some $r\geq 1$, then $\Xi^{N,t}_{m, q}=\emptyset$ for all $q\geq r$. Let
\[
\Gamma_{m}^{N, t}=\bigcup_{r=0}^{N-1}\Xi^{N, t}_{m , r},
\]
then $\xi_t^N(m)$ is determined by values of $\{\xi_0^N(\varphi)\}_{\varphi\in \Gamma_{m}^{N, t}}$, $\{K_\varphi^N\}_{\varphi\in \Gamma_{m}^{N, t}}$ and $\{U^N_{(\varphi, \varrho)}:~\varphi\in \Gamma_{m}^{N, t}, \varrho\neq \varphi\}$ since it takes an infected vertex out of $\Gamma_{m}^{N, t}$ longer than $t$ to infect $m$ through an infection path.

For given four different $i,j,k,l\in \{1,2,\ldots, N\}$, let $\{\Gamma_{m}^{N, t}:~m=i,j,k,l\}$ be defined as above. Now we define the copy
$\hat{\Gamma}_m^{N,t}$ of $\Gamma_m^{N, t}$ for $m=j,k,l$. For $m=j,k,l$, let
\[
\hat{\Gamma}_{m}^{N, t}=\bigcup_{r=0}^{N-1}\hat{\Xi}^{N, t}_{m , r},
\]
where $\{\hat{\Xi}^{N, t}_{m, r}\}_{r=0}^{N-1}$ is defined by induction as follows:

1) $\hat{\Xi}^{N, t}_{m, 0}=\{m\}$.

2) After $\hat{\Xi}^{N,t}_{m, r}$ are defined for $0\leq r\leq q<N-1$, $\hat{\Xi}^{N, t}_{m, q+1}$ is defined as the set of $\varphi\in \{1,2,\ldots, N\}\setminus\left(\bigcup_{r=0}^q\hat{\Xi}^{N,t}_{m, r}\right)$ such that there exists self-avoiding path $\varphi=m_{q+1}, m_q, \ldots, m_1, m_0=m$ such that $m_r\in \bigcup_{l=1}^r\hat{\Xi}^{N, t}_{m,l}$ for $1\leq r\leq q$ and
\[
\sum_{r=0}^{q}\tilde{U}^{N, m}_{(m_r, m_{r+1})}\leq t,
\]
where for $m=j,k,l$ and all $1\leq \varphi\neq \varrho\leq N$,
\[
\tilde{U}^{N, m}_{(\varphi, \varrho)}=
\begin{cases}
U^N_{(\varphi, \varrho)} & \text{~if~} m=j\text{~and~}\varphi\not\in \Gamma_{i}^{N, t}, \\
U^N_{(\varphi, \varrho), 2} & \text{~if~} m=j\text{~and~}\varphi \in \Gamma_{i}^{N, t}, \\
U^N_{(\varphi, \varrho)} & \text{~if~} m=k\text{~and~}\varphi \not\in \Gamma_{i}^{N, t}\bigcup \hat{\Gamma}_{j}^{N, t}, \\
U^N_{(\varphi, \varrho), 3} & \text{~if~} m=k\text{~and~}\varphi \in \Gamma_{i}^{N, t}\bigcup \hat{\Gamma}_{j}^{N, t},\\
U^N_{(\varphi, \varrho)} & \text{~if~} m=l\text{~and~}\varphi\not\in  \Gamma_{i}^{N, t}\bigcup \hat{\Gamma}_{j}^{N, t}\bigcup \hat{\Gamma}_{k}^{N, t}, \\
U^N_{(\varphi, \varrho), 4} & \text{~if~} m=l\text{~and~}\varphi\in  \Gamma_{i}^{N, t}\bigcup \hat{\Gamma}_{j}^{N, t}\bigcup \hat{\Gamma}_{k}^{N, t}.
\end{cases}
\]
We claim that $\hat{\Gamma}_{m}^{N, t}$ and $\Gamma_m^{N, t}$ have the same distribution for $m=j,k,l$ and
\[
\Gamma_i^{N,t}, \hat{\Gamma}_j^{N, t}, \hat{\Gamma}_k^{N, t}, \hat{\Gamma}_l^{N, t}
\]
are independent. This claim is easy to check since $\Gamma_i^{N,t}, \hat{\Gamma}_j^{N, t}, \hat{\Gamma}_k^{N, t}, \hat{\Gamma}_l^{N, t}$ can be obtained from  a breadth-first search (BFS). In detail, to obtain $\Xi^{N, t}_{i, 1}$, we check values of $U^N_{(i, \varphi)}$  for all $\varphi\neq i$ and then $\Xi^{N, t}_{i, 1}$ is the set of $\varphi$ satisfying $U^N_{(i, \varphi)}\leq t$. Then to obtain $\Xi_{i, 2}^{N, t}$, we check values of $U^N_{(\varphi, \varrho)}$ for all $\varphi\in \Xi^{N, t}_{i, 1}, \varphi\neq \varrho$ and then $\Xi^{N, t}_{i, 2}$ is the set of $\varrho\not\in \{i\}\bigcup\Xi^{N, t}_{i, 1}$ such that there exists $\varphi\in \Xi^{N, t}_{i, 1}$ such that $U^N_{(i, \varphi)}+U^N_{(\varphi, \varrho)}\leq t$. For $r\geq 3$, we obtain $\Xi^{N, t}_{i, r}$ by repeating the above procedure. After $\Gamma_i^{N, t}$ is obtained, we can obtain $\hat{\Gamma}_j^{N, t}$ according to the same BFS procedure. However, in the second BFS if some $\varphi\in \Gamma_i^{N, t}$ is included in $\bigcup_{l=1}^r\hat{\Xi}^{N, t}_{j, l}$, then in the $r+1$th step, we check the value of $U^N_{(\varphi, \varrho), 2}$ for $\varrho\neq \varphi$ instead of $U^N_{(\varphi, \varrho)}$. According to this rule, the information of an exponential infection time is utilized in at most one of the definition of $\Gamma_i^{N, t}$ and the definition of $\hat{\Gamma}_j^{N, t}$ and hence $\hat{\Gamma}_j^{N, t}$ is independent of $\Gamma_i^{N, t}$. Similar analysis holds for $\hat{\Gamma}_k^{N, t}$ and $\hat{\Gamma}_l^{N, t}$.

For $m=j,k,l$, let $\hat{\xi}_t^N(m)$ be determined by values of $\{\tilde{\xi}_0^{N, m}(\varphi)\}_{\varphi\in \hat{\Gamma}_{m}^{N, t}}$, $\{\tilde{K}_\varphi^{N, m}\}_{\varphi\in \hat{\Gamma}_{m}^{N, t}}$ and $\{\tilde{U}^{N,m}_{(\varphi, \varrho)}:~\varphi\in \hat{\Gamma}_{m}^{N, t}, \varrho\neq \varphi\}$ in the same way as $\xi_t^N(m)$ is determined by values of $\{\xi_0^N(\varphi)\}_{\varphi\in \Gamma_{m}^{N, t}}$, $\{K_\varphi^N\}_{\varphi\in \Gamma_{m}^{N, t}}$ and $\{U^N_{(\varphi, \varrho)}:~\varphi\in \Gamma_{m}^{N, t}, \varrho\neq \varphi\}$, where
\[
\tilde{\xi}_0^{N, m}(\varphi)=
\begin{cases}
\xi_0^N(\varphi) & \text{~if~} m=j\text{~and~}\varphi\not\in \Gamma_{i}^{N, t}, \\
\xi_{0,2}^N(\varphi) & \text{~if~} m=j\text{~and~}\varphi \in \Gamma_{i}^{N, t}, \\
\xi_0^N(\varphi) & \text{~if~} m=k\text{~and~}\varphi \not\in \Gamma_{i}^{N, t}\bigcup \hat{\Gamma}_{j}^{N, t}, \\
\xi_{0,3}^N(\varphi) & \text{~if~} m=k\text{~and~}\varphi \in \Gamma_{i}^{N, t}\bigcup \hat{\Gamma}_{j}^{N, t},\\
\xi_0^N(\varphi) & \text{~if~} m=l\text{~and~}\varphi\not\in  \Gamma_{i}^{N, t}\bigcup \hat{\Gamma}_{j}^{N, t}\bigcup \hat{\Gamma}_{k}^{N, t}, \\
\xi_{0,4}^N(\varphi) & \text{~if~} m=l\text{~and~}\varphi\in  \Gamma_{i}^{N, t}\bigcup \hat{\Gamma}_{j}^{N, t}\bigcup \hat{\Gamma}_{k}^{N, t}
\end{cases}
\]
and
\[
\tilde{K}_\varphi^{N, m}=
\begin{cases}
K_\varphi^N & \text{~if~} m=j\text{~and~}\varphi\not\in \Gamma_{i}^{N, t}, \\
K_{\varphi,2}^N & \text{~if~} m=j\text{~and~}\varphi \in \Gamma_{i}^{N, t}, \\
K_\varphi^N & \text{~if~} m=k\text{~and~}\varphi \not\in \Gamma_{i}^{N, t}\bigcup \hat{\Gamma}_{j}^{N, t}, \\
K_{\varphi,3}^N & \text{~if~} m=k\text{~and~}\varphi \in \Gamma_{i}^{N, t}\bigcup \hat{\Gamma}_{j}^{N, t},\\
K_\varphi^N & \text{~if~} m=l\text{~and~}\varphi\not\in  \Gamma_{i}^{N, t}\bigcup \hat{\Gamma}_{j}^{N, t}\bigcup \hat{\Gamma}_{k}^{N, t}, \\
K_{\varphi,4}^N & \text{~if~} m=l\text{~and~}\varphi\in  \Gamma_{i}^{N, t}\bigcup \hat{\Gamma}_{j}^{N, t}\bigcup \hat{\Gamma}_{k}^{N, t}.
\end{cases}
\]
Since the information of each exponential recovery time or initial state is utilized in at most one of respective definitions of $\xi_t^N(i), \hat{\xi}_t^N(j), \hat{\xi}_t^N(k), \hat{\xi}_t^N(l)$, these four random variables are independent. Now we can give our proof of Lemma \ref{lemma 3.1 replacement}.

\proof[Proof of Lemma \ref{lemma 3.1 replacement}]

Throughout this proof we assume that $i,j,k,l$ are four given different vertices in $\{1,2,\ldots, N\}$ and $0\leq t\leq T$. We only give the proof of Equation \eqref{equ 3.1 replacement lemma partThree} since the other two follows from similar analysis.

For any $B\subset \{1,2,\ldots, N\}$ and $m\not\in B$, we denote by $\Xi^{N,t,B}_{m, r}$ the result of the $r$th step of the BFS as that giving $\Gamma_m^{N,t}$ except that vertices in $B$ are not investigated, i.e., we have:

1) $\Xi^{N, t, B}_{m, 0}=\{m\}$.

2) After $\Xi^{N,t, B}_{m, r}$ are defined for $0\leq r\leq q<N-1$, $\Xi^{N, t, B}_{m, q+1}$ is defined as the set of $\varphi\in \left(\{1,2,\ldots, N\}\setminus B\right)\setminus\left(\bigcup_{r=0}^q\Xi^{N, t, B}_{m, r}\right)$ such that there exists self-avoiding path
\[
\varphi=m_{q+1}, m_q, \ldots, m_1, m_0=m
\]
such that $m_r\in \bigcup_{l=1}^r\Xi^{N, t, B}_{m,l}$ for $1\leq r\leq q$ and
\[
\sum_{r=0}^{q}U^N_{(m_r, m_{r+1})}\leq t.
\]
We further define $\Gamma_m^{N, t, B}=\bigcup_{r=0}^{N-1}\Xi^{N,t,B}_{m, r}$. For $m=j,k,l$ and $B\not\ni m$, $\hat{\Xi}_{m,r}^{N, t, B}$ and $\hat{\Gamma}_m^{N,t,B}$ are defined according to $\{\tilde{U}^{N,m}_{(\varphi, \varrho)}\}_{1\leq \varphi\neq \varrho\leq N}$ in the same way as $\Xi^{N,t,B}_{m, r}$ and $\Gamma_m^{N, t, B}$ are defined according to $\{U^N_{(\varphi, \varrho)}\}_{1\leq \varphi\neq \varrho\leq N}$ . For simplicity, we define
\[
B_1^N=\Gamma_i^{N, t, \{j, k, l\}}, B_2^N=\hat{\Gamma}_j^{N, t, B_1^N\bigcup\{k,l\}}, B_3^N=\hat{\Gamma}_k^{N, t, \bigcup_{m=1}^2B_m^N\bigcup \{l\}}
\]
and
\[
B_4^N=\hat{\Gamma}_l^{N, t, \bigcup_{m=1}^3B_m^N}.
\]
We denote by $\Omega_{N, 1}$ the event that $U^N_{(\varphi, \varrho)}>T$ for any $(\varphi, \varrho)\in B_{m_1}^N\times B_{m_2}^N$ for all $1\leq m_2<m_1\leq 4$, by $\Omega_{N, 2}$ the event that $U^N_{(\varphi, j)}>T$ for any $\varphi\in B_1^N$, by $\Omega_{N, 3}$ the event that $U^N_{(\varphi, k)}>T$ for any $\varphi\in B_1^N\bigcup B_2^N$, by $\Omega_{N, 4}$ the event that $U^N_{(\varphi, l)}>T$ for any $\varphi\in B_1^N\bigcup B_2^N\bigcup B_3^N$ and by $\Omega_N$ the intersection of $\Omega_{N, m}$ for $m=1,2,3,4$. On the event $\Omega_N$,
\[
B_1^N=\Gamma_i^{N, t}, B_2^N=\hat{\Gamma}_j^{N, t}=\Gamma_j^{N, t}, B_3^N=\hat{\Gamma}_k^{N, t}=\Gamma_k^{N, t}
\]
and
\[
B_4^N=\hat{\Gamma}_l^{N, t}=\Gamma_l^{N, t}.
\]
Consequently, on the event $\Omega_N$, $\Gamma_{m_1}^{N,t}\bigcap \Gamma_{m_2}^{N, t}=\emptyset$ for any different $m_1, m_2\in \{i,j,k,l\}$ and $\hat{\xi}_t^N(m)=\xi_t^N(m)$ for $m=j,k,l$. As a result, for $(r_1, r_2, r_3, r_4)=(1,1,1,1)$,
\begin{align}\label{equ 3.4}
&\left|\mathbb{E}(H_{t, r_1}^N(i)H_{t, r_2}^N(j)H_{t, r_3}^N(k)H_{t, r_4}^N(l))-\mathbb{E}H_{t, r_1}^N(i)\mathbb{E}H_{t, r_2}^N(j)\mathbb{E}H_{t, r_3}^N(k)\mathbb{E}H_{t, r_4}^N(l)\right| \notag\\
&=\left|P\left(\xi_t^N(m)=1\text{~for~}m=i,j,k,l\right)-\prod_{m=i,j,k,l}P\left(\xi_t^N(m)=1\right)\right| \notag\\
&=\left|P\left(\xi_t^N(m)=1\text{~for~}m=i,j,k,l\right)-P(\xi_t^N(i)=1)\prod_{m=j,k,l}P\left(\hat{\xi}_t^N(m)=1\right)\right|\notag\\
&=\left|P\left(\xi_t^N(m)=1\text{~for~}m=i,j,k,l\right)-P\left(\xi_t(i)=1, \hat{\xi}_t^N(m)=1\text{~for~}m=j,k,l\right)\right|\notag\\
&\leq \Bigg|P\left(\xi_t^N(m)=1\text{~for~}m=i,j,k,l, \Omega_N^c\right)\notag\\
&\text{\quad\quad\quad}-P\left(\xi_t(i)=1, \hat{\xi}_t^N(m)=1\text{~for~}m=j,k,l, \Omega_N^c\right)\Bigg|\notag\\
&\leq 2P(\Omega_N^c).
\end{align}
Equation \eqref{equ 3.4} holds for other $(r_1, r_2, r_3, r_4)\in \{1, 0\}^4$ according to a similar analysis. For any $\varrho\neq\varphi$,
\[
P\left(U^N_{(\varrho, \varphi)}\leq T\right)=1-e^{-\frac{1}{N}\lambda\left(\frac{i}{N}, \frac{j}{N}\right)T}\leq \frac{\|\lambda\|_\infty T}{N},
\]
where $\|\lambda\|_\infty=\sup_{0\leq u,v\leq 1}|\lambda(u, v)|$. Note that information of
\[
\{U^N_{(\varphi, \varrho)}:~(\varphi, \varrho)\in B_{m_1}^N\times B_{m_2}^N\text{~for some~}m_2<m_1\}, \{U^N_{(\varphi, j)}\}_{\varphi\in B_1^N},
\]
$\{U^N_{(\varphi, k)}\}_{\varphi\in B_1^N\bigcup B_2^N}$ and $\{U^N_{(\varphi, l)}\}_{\varphi\in B_1^N\bigcup B_2^N\bigcup B_3^N}$ are not utilized in the definition of $B_1^N, B_2^N, B_3^N, B_4^N$, as a result,
\[
P(\Omega_N^c)\leq \frac{\|\lambda\|_\infty T}{N}\mathbb{E}\left(3|B_1^N|+2|B_2^N|+|B_3^N|+\sum_{m_2<m_1}|B_{m_2}^N||B_{m_1}^N|\right),
\]
where $|A|$ is the cardinality of the set $A$. For any self-avoiding path $m_0, m_1,\ldots, m_q$ on $\{1,2,\ldots,N\}$ with length $q$ and $\tau>0$,
\begin{align*}
P\left(\sum_{r=0}^{q-1}U^N_{(m_r, m_{r+1})}\leq t\right)\leq e^{\tau t}\prod_{r=0}^{q-1}\mathbb{E}e^{-\tau U^N_{(m_r, m_{r+1})}}\leq e^{\tau t}\left(\frac{\|\lambda\|_\infty}{N\tau+\|\lambda\|_\infty}\right)^q.
\end{align*}
Let $\tau=2\|\lambda\|_\infty$, then for any $m\in \{1,2,\ldots, N\}$,
\[
\mathbb{E}|\Gamma_m^{N, t}|\leq \sum_{q=0}^{+\infty}N^qe^{2\|\lambda\|_\infty T}\left(\frac{1}{2N}\right)^q=2e^{2\|\lambda\|_\infty T}.
\]
For $(m_2, m_1)=(1,2)$,
\begin{align*}
\mathbb{E}\left(|B_1^N||B_2^N|\right)&\leq \mathbb{E}\left(|\Gamma_i^{N,t}||\hat{\Gamma}^{N, t}_j|\right)=\mathbb{E}|\Gamma_i^{N,t}|\mathbb{E}|\hat{\Gamma}_j^{N, t}|\leq 4e^{4\|\lambda\|_\infty T}.
\end{align*}
The above inequality holds for other $m_2<m_1$ according to similar analysis. As a result,
\begin{equation}\label{equ 3.5}
P(\Omega_N^c)\leq \frac{\|\lambda\|_\infty T}{N}\left(12e^{2\|\lambda\|_\infty T}+24e^{4\|\lambda\|_\infty T}\right).
\end{equation}
Equation \eqref{equ 3.1 replacement lemma partThree} follows from Equations \eqref{equ 3.4} and \eqref{equ 3.5}.

\qed

\section{Proof of Theorem \ref{theorem 2.1 hydrodynamic limit}}\label{section four}

In this section we prove Theorem \ref{theorem 2.1 hydrodynamic limit}. It is easy to check that Theorem \ref{theorem 2.1 hydrodynamic limit} follows from the following two lemmas.

\begin{lemma}\label{lemma 4.1}
For any $f\in C[0, 1]$,
\[
\lim_{N\rightarrow+\infty}{\rm Var}\left(\mu_t^N(f)\right)=\lim_{N\rightarrow+\infty}{\rm Var}\left(\theta_t^N(f)\right)=0.
\]
\end{lemma}

\begin{lemma}\label{lemma 4.2}
For any $f\in C[0, 1]$,
\[
\lim_{N\rightarrow+\infty}\mathbb{E}\left(\mu_t^N(f)\right)=\int_0^1\rho_1(t,u)f(u)du
\]
and
\[
\lim_{N\rightarrow+\infty}\mathbb{E}\left(\theta_t^N(f)\right)=\int_0^1\rho_0(t, u)f(u)du.
\]
\end{lemma}

\proof[Proof of Theorem \ref{theorem 2.1 hydrodynamic limit}]

According to the fact that $(c+d)^2\leq 2c^2+2d^2$ for $c,d\in \mathbb{R}$,
\begin{align*}
&\left(\mu_t^N(f)-\int_0^1\rho_1(t,u)f(u)du\right)^2\\
&\leq 2\left(\mu_t^N(f)-\mathbb{E}\left(\mu_t^N(f)\right)\right)^2+2\left(\mathbb{E}\left(\mu_t^N(f)\right)-\int_0^1\rho_1(t,u)f(u)du\right)^2.
\end{align*}
Then, by Lemmas \ref{lemma 4.1} and \ref{lemma 4.2}, $\mu_t^N(f)$ converges to $\int_0^1\rho_1(t,u)f(u)du$ in $L^2$ and hence in probability as $N\rightarrow+\infty$. Similar analysis holds for $\theta_t^N(f)$.

\qed

The remainder of this section is devoted to proofs of Lemmas \ref{lemma 4.1} and \ref{lemma 4.2}.

\proof[Proof of Lemma \ref{lemma 4.1}]

Since $T$ can be chosen arbitrarily, we assume that $0<t\leq T$ without loss of generality. By Lemma \ref{lemma 3.1 replacement},

\begin{align*}
&{\rm Var}\left(\mu_t^N(f)\right)=\frac{1}{N^2}\sum_{i=1}^N\sum_{j=1}^N{\rm Cov}\left(I_t^N(i), I_t^N(j)\right)f\left(\frac{i}{N}\right)f\left(\frac{j}{N}\right)\\
&\leq \frac{\|f\|_\infty^2}{N^2}\sum_{i=1}^N\sum_{j=1}^N\left|\mathbb{E}(H_{t, 1}^N(i)H_{t, 1}^N(j))-\mathbb{E}H_{t, 1}^N(i)\mathbb{E}H_{t, 1}^N(j)\right|\\
&=\frac{\|f\|_\infty^2}{N^2}\sum_{i=1}^N\sum_{j\neq i}\left|\mathbb{E}(H_{t, 1}^N(i)H_{t, 1}^N(j))-\mathbb{E}H_{t, 1}^N(i)\mathbb{E}H_{t, 1}^N(j)\right|\\
&\text{\quad\quad}+\frac{\|f\|_\infty^2}{N^2}\sum_{i=1}^N\left|\mathbb{E}(H_{t, 1}^N(i)H_{t, 1}^N(i))-\mathbb{E}H_{t, 1}^N(i)\mathbb{E}H_{t, 1}^N(i)\right|\\
&\leq \frac{2N\|f\|_\infty^2}{N^2}+\frac{\|f\|_\infty^2N(N-1)C_1}{N^3}=\|f\|_\infty^2\left(\frac{(N-1)C_1}{N^2}+\frac{2}{N}\right),
\end{align*}
where $\|f\|_\infty=\sup_{0\leq u\leq 1}|f(u)|$. As a result,
\[
\lim_{N\rightarrow+\infty}{\rm Var}\left(\mu_t^N(f)\right)=0.
\]
The second part of Lemma \ref{lemma 4.1} that $\lim_{N\rightarrow+\infty}{\rm Var}\left(\theta_t^N(f)\right)=0$ follows from a similar analysis.

\qed

\proof[Proof of Lemma \ref{lemma 4.2}]

We still assume that $0<t\leq T$ without loss of generality. For $N\geq 1$, $1\leq i\leq N$ and any $0\leq t\leq T$, let
\[
\rho_1^N(t, u)=P(I_t^N(i)=1)\text{~and~}\rho_0^N(t, u)=P(S_t^N(i)=1)
\]
when $\frac{i-1}{N}<u\leq \frac{i}{N}$. Then, to prove Lemma \ref{lemma 4.2} we only need to show that
\begin{equation}\label{equ 4.1}
\lim_{N\rightarrow+\infty}\sup_{0\leq t\leq T}\sup_{0<u\leq 1}\left|\rho_m^N(t, u)-\rho_m(t, u)\right|=0
\end{equation}
for $m=0, 1$. According to the generator $\mathcal{L}^N$ of $\{\xi_t^N\}_{t\geq 0}$ given in Section \ref{section one} and Chapman-Kolmogorov Equation, for any $1\leq i\leq N$ and $\frac{i-1}{N}<u\leq \frac{i}{N}$,
\begin{align*}
\frac{d}{dt}\rho_1^N(t, u)&=\frac{d}{dt}P(I_t^N(i)=1) \\
&=-\psi\left(\frac{i}{N}\right)P(I_t^N(i)=1)+\frac{1}{N}\sum_{j\neq i}\lambda\left(\frac{i}{N}, \frac{j}{N}\right)P(S_t^N(i)=1, I_t^N(j)=1)
\end{align*}
and
\begin{align*}
\frac{d}{dt}\rho_0^N(t, u)&=\frac{d}{dt}P(S_t^N(i)=1) \\
&=-\frac{1}{N}\sum_{j\neq i}\lambda\left(\frac{i}{N}, \frac{j}{N}\right)P(S_t^N(i)=1, I_t^N(j)=1).
\end{align*}
For all $N\geq 1$, let
\[
\Delta_\lambda^N=\sup\left\{|\lambda(u_1, v_1)-\lambda(u_2, v_2)|:~|u_1-u_2|, |v_1-v_2|\leq N^{-1}\right\},
\]
\[
\Delta_\psi^N=\sup\left\{|\psi(u_1)-\psi(u_2)|:~|u_1-u_2|\leq N^{-1}\right\}
\]
and
\[
\Delta_\phi^N=\sup\left\{|\phi(u_1)-\phi(u_2)|:~|u_1-u_2|\leq N^{-1}\right\}.
\]
By Lemma \ref{lemma 3.1 replacement}, $\left|P(S_t^N(i)=1, I_t^N(j)=1)-\rho_0^N\left(t, \frac{i}{N}\right)\rho_t^N\left(t, \frac{j}{N}\right)\right|\leq \frac{C_1}{N}$ for $i\neq j$ and hence
\begin{align*}
&\left|\frac{1}{N}\sum_{j\neq i}\lambda\left(\frac{i}{N}, \frac{j}{N}\right)P(S_t^N(i)=1, I_t^N(j)=1)
-\rho_0^N(t,u)\int_0^1\lambda(u, v)\rho_1^N(t,v)dv\right|\\
&\leq \Delta_\lambda^N+\frac{C_1\|\lambda_\infty\|}{N}+\frac{\|\lambda\|_\infty}{N}
\end{align*}
and $|\psi\left(\frac{i}{N}\right)P(I_t^N(i)=1)-\psi(u)\rho_1^N(t,u)|\leq \Delta_\psi^N$ for $\frac{i-1}{N}<u\leq \frac{i}{N}$. For all $N\geq 1$, let
\[
M_t^N=\sup_{0<u\leq 1}\left\{|\rho_0^N(t,u)-\rho_0(t,u)|+|\rho_1^N(t,u)-\rho_1(t, u)|\right\},
\]
then
\[
|\rho_0^N(t,u)\lambda(u, v)\rho_1^N(t,v)-\rho_0(t,u)\lambda(u, v)\rho_1(t,v)|\leq 2\|\lambda\|_\infty M_t^N
\]
for any $0<u,v\leq 1$. According to Assumption (A), $M_0^N\leq 2\Delta_\phi^N$. In conclusion, by writing Equation \eqref{equ 2.1 nonliear ODE} and above two Chapman-Kolmogorov Equations as their integration forms, we have
\[
M_t^N\leq \left(2\Delta_\phi^N+T(\Delta_\psi^N+2\Delta_\lambda^N+\frac{2C_1\|\lambda_\infty\|}{N}+\frac{2\|\lambda\|_\infty}{N})\right)+4\|\lambda\|_\infty\int_0^tM_s^Nds
\]
for $0\leq t\leq T$. As a result, according to Gronwall's inequality,
\[
M_t^N\leq \left(2\Delta_\phi^N+T(\Delta_\psi^N+2\Delta_\lambda^N+\frac{2C_1\|\lambda_\infty\|}{N}+\frac{2\|\lambda\|_\infty}{N})\right)e^{4\|\lambda\|_\infty T}
\]
for $0\leq t\leq T$. Consequently, $\lim_{N\rightarrow+\infty}\sup_{0\leq t\leq T}M_t^N=0$ and hence Equation \eqref{equ 4.1} holds.

\qed

\section{Proof of Theorem \ref{theorem 2.2 fluctuation}}\label{section five}

In this section we prove Theorem \ref{theorem 2.2 fluctuation}. As a preliminary, we first need the tightness of $\{\left(\eta_t^N, \beta_t^N\right):~0\leq t\leq T\}_{N\geq 1}$.

\begin{lemma}\label{lemma 5.1}
Under Assumption (A), $\{\left(\eta_t^N, \beta_t^N\right):~0\leq t\leq T\}_{N\geq 1}$ are tight.
\end{lemma}

\proof

According to Aldous' criteria, the tightness of $\{\left(\eta_t^N, \beta_t^N\right):~0\leq t\leq T\}_{N\geq 1}$ follows from the following four equations.

1) For all $0\leq t\leq T$ and $f\in C[0, 1]$,
\begin{equation}\label{equ 5.1}
\lim_{M\rightarrow+\infty}\limsup_{N\rightarrow+\infty}P\left(\left|\beta_t^N(f)\right|\geq M\right)=0.
\end{equation}

2) For all $\epsilon>0$ and $f\in C[0, 1]$,
\begin{equation}\label{equ 5.2}
\lim_{\delta\rightarrow0}\limsup_{N\rightarrow+\infty}\sup_{\varsigma\in \mathcal{T}, s<\delta}P\left(\left|\beta^N_{\varsigma+s}(f)-\beta_\varsigma^N(f)\right|>\epsilon\right)=0,
\end{equation}
where $\mathcal{T}$ is the set of stopping times of $\{\xi_t^N\}_{t\geq 0}$ bounded by $T$.

3) For all $0\leq t\leq T$ and $f\in C[0, 1]$,
\begin{equation*}
\lim_{M\rightarrow+\infty}\limsup_{N\rightarrow+\infty}P\left(\left|\eta_t^N(f)\right|\geq M\right)=0.
\end{equation*}

4) For all $\epsilon>0$ and $f\in C[0, 1]$,
\begin{equation*}
\lim_{\delta\rightarrow0}\limsup_{N\rightarrow+\infty}\sup_{\varsigma\in \mathcal{T}, s<\delta}P\left(\left|\eta^N_{\varsigma+s}(f)-\eta_\varsigma^N(f)\right|>\epsilon\right)=0.
\end{equation*}

Here we only check Equations \eqref{equ 5.1} and \eqref{equ 5.2} since the other two follows from a similar analysis. According to Dynkin's martingale formula,
\[
\beta_t^N(f)=\beta_0^N(f)+\int_0^t\left(\partial_s+\mathcal{L}^N\right)\beta_s^N(f)ds+\mathcal{M}_{t,\beta}^N(f),
\]
where $\{\mathcal{M}_{t,\beta}^N(f)\}_{t\geq 0}$ is a martingale with mean zero and quadratic variation process $\langle\mathcal{M}^N_\beta(f)\rangle_t$ given by
\[
\langle\mathcal{M}^N_\beta(f)\rangle_t=\int_0^t\mathcal{L}^N\left(\left(\beta_s^N(f)\right)^2\right)-2\beta_s^N(f)\mathcal{L}^N\beta_s^N(f)ds.
\]
As a result, to prove Equations \eqref{equ 5.1} and \eqref{equ 5.2}, we only need to check the following four equations.

1) For all $0\leq t\leq T$ and $f\in C[0, 1]$,
\begin{equation}\label{equ 5.3}
\lim_{M\rightarrow+\infty}\limsup_{N\rightarrow+\infty}P\left(\left|\int_0^t\left(\partial_s+\mathcal{L}^N\right)\beta_s^N(f)ds\right|\geq M\right)=0.
\end{equation}

2) For all $\epsilon>0$ and $f\in C[0, 1]$,
\begin{equation}\label{equ 5.4}
\lim_{\delta\rightarrow0}\limsup_{N\rightarrow+\infty}\sup_{\varsigma\in \mathcal{T}, s<\delta}P\left(\left|\int_\varsigma^{\varsigma+s} \left(\partial_u+\mathcal{L}^N\right)\beta_u^N(f)du\right|>\epsilon\right)=0.
\end{equation}

3) For all $0\leq t\leq T$ and $f\in C[0, 1]$,
\begin{equation}\label{equ 5.5}
\lim_{M\rightarrow+\infty}\limsup_{N\rightarrow+\infty}P\left(\left|\mathcal{M}_{t, \beta}^N(f)\right|\geq M\right)=0.
\end{equation}

4) For all $\epsilon>0$ and $f\in C[0, 1]$,
\begin{equation}\label{equ 5.6}
\lim_{\delta\rightarrow0}\limsup_{N\rightarrow+\infty}\sup_{\varsigma\in \mathcal{T}, s<\delta}P\left(\left|\mathcal{M}_{\varsigma+s, \beta}^N(f)-\mathcal{M}_{\varsigma, \beta}^N(f)\right|>\epsilon\right)=0.
\end{equation}

For Equation \eqref{equ 5.3}, by direct calculation,
\begin{align*}
&\int_0^t\left(\partial_s+\mathcal{L}^N\right)\beta_s^N(f)ds \\
&=\int_0^t\left(\frac{1}{N^{1.5}}\sum_{i=1}^N\sum_{j\neq i}\lambda\left(\frac{i}{N}, \frac{j}{N}\right)\left(S_t^N(i)I_t^N(j)-\mathbb{E}\left(S_t^N(i)I_t^N(j)\right)\right)f\left(\frac{i}{N}\right)\right)ds.
\end{align*}

According to Lemma \ref{lemma 3.1 replacement}, for any different $i,j,k,l$,
\begin{equation}\label{equ 5.7}
|{\rm Cov}\left(S_t^N(i)I_t^N(j), S_t^N(k)I_t^N(l)\right)|\leq \frac{C_3+2C_1}{N}+\frac{C_1^2}{N^2}\leq \frac{C_5}{N},
\end{equation}
where $C_5=C_3+2C_1+C_1^2$. As a result, by Cauchy-Schwarz inequality,
\begin{align}\label{equ 5.8}
&\mathbb{E}\left(\left(\int_0^t\left(\partial_s+\mathcal{L}^N\right)\beta_s^N(f)ds\right)^2\right)\notag\\
&\leq T^2\|\lambda\|_\infty^2\|f\|_\infty^2\left(\frac{C_5}{N}\frac{N(N-1)(N-2)(N-3)}{N^3}+\frac{6N^3}{N^3}\right) \notag\\
&\leq T^2\|\lambda\|_\infty^2\|f\|_\infty^2\left(C_5+6\right).
\end{align}
Equation \eqref{equ 5.3} follows from Equation \eqref{equ 5.8} and Markov's inequality.

\quad

For Equation \eqref{equ 5.4}, by Cauchy-Schwarz inequality and Equation \eqref{equ 5.7},
\begin{align}\label{equ 5.9}
{\mathbb{E}}\left(\left|\int_\varsigma^{\varsigma+s} \left(\partial_u+\mathcal{L}^N\right)\beta_u^N(f)du\right|^2\right)
&\leq \delta\int_0^{T+\delta}\mathbb{E}\left(\left(\left(\partial_u+\mathcal{L}^N\right)\beta_u^N(f)\right)^2\right)du \notag\\
&\leq \delta(T+\delta)\|\lambda\|_\infty^2\|f\|_\infty^2(C_5+6).
\end{align}
Equation \eqref{equ 5.4} follows from Equation \eqref{equ 5.9} and Markov's inequality.

\quad

For Equation \eqref{equ 5.5}, by direct calculation,
\begin{align}\label{equ 5.10}
&\mathcal{L}^N\left(\left(\beta_s^N(f)\right)^2\right)-2\beta_s^N(f)\mathcal{L}^N\beta_s^N(f) \notag\\
&=\frac{1}{N^2}\sum_{i=1}^N\sum_{j\neq i}\lambda\left(\frac{i}{N}, \frac{j}{N}\right)S_s^N(i)I_s^N(j)f^2\left(\frac{i}{N}\right)\leq \|\lambda\|_\infty\|f\|_\infty^2.
\end{align}
Hence,
\begin{equation}\label{equ 5.11}
\mathbb{E}\left(\left|\mathcal{M}_{t, \beta}^N(f)\right|^2\right)=\mathbb{E}\langle\mathcal{M}^N_\beta(f)\rangle_t\leq T\|\lambda\|_\infty\|f\|^2_\infty.
\end{equation}
Equation \eqref{equ 5.5} follows from Equation \eqref{equ 5.11} and Markov's inequality.

\quad

For Equation \eqref{equ 5.6}, by Equation \eqref{equ 5.10},
\begin{align}\label{equ 5.12}
\mathbb{E}\left(\left|\mathcal{M}_{\varsigma+s, \beta}^N(f)-\mathcal{M}_{\varsigma, \beta}^N(f)\right|^2\right)
&=\mathbb{E}\left(\langle\mathcal{M}^N_\beta(f)\rangle_{\varsigma+s}-\langle\mathcal{M}^N_\beta(f)\rangle_{\varsigma}\right) \notag\\
&\leq \delta \|\lambda\|_\infty\|f\|^2_\infty.
\end{align}
Equation \eqref{equ 5.6} follows from Equation \eqref{equ 5.12} and Markov's inequality.

In conclusion, Equations \eqref{equ 5.3}-\eqref{equ 5.6} all hold and hence Equations \eqref{equ 5.1} and \eqref{equ 5.2} hold.

\qed

For later use, we need the following lemma.

\begin{lemma}\label{lemma 5.2}
For any $h\in C\left([0, 1]\times[0, 1]\right)$ and $0\leq t\leq T$,
\[
\lim_{N\rightarrow+\infty}\int_0^t\frac{1}{N^{1.5}}\sum_{i=1}^N\sum_{j\neq i}\left(S_s^N(i)-\mathbb{E}S_s^N(i)\right)\left(I_s^N(j)-\mathbb{E}I_s^N(j)\right)h\left(\frac{i}{N}, \frac{j}{N}\right)ds=0
\]
in probability.
\end{lemma}

\proof

Since $S_t^N(i), I_t^N(i)\leq 1$, we only need to show that
\begin{equation}\label{equ 5.13}
\lim_{N\rightarrow+\infty}\int_0^t\frac{1}{N^{1.5}}\sum_{i=1}^N\sum_{j=1}^N
\left(S_s^N(i)-\mathbb{E}S_s^N(i)\right)\left(I_s^N(j)-\mathbb{E}I_s^N(j)\right)h\left(\frac{i}{N}, \frac{j}{N}\right)ds=0
\end{equation}
in probability. Let $\mathcal{H}$ be the set of $h\in C\left([0, 1]\times[0, 1]\right)$ that $h(u,v)=h_1(u)h_2(v)$ for all $u,v\in [0, 1]$ for some $h_1, h_2\in C[0, 1]$. We first show that Equation \eqref{equ 5.13} holds for all $h\in\mathcal{H}$ and hence for all $h\in{\rm span}(\mathcal{H})$. For any $h\in \mathcal{H}$,
\begin{align*}
&\mathbb{E}\left|\frac{1}{N^{1.5}}\sum_{i=1}^N\sum_{j=1}^N
\left(S_s^N(i)-\mathbb{E}S_s^N(i)\right)\left(I_s^N(j)-\mathbb{E}I_s^N(j)\right)h\left(\frac{i}{N}, \frac{j}{N}\right)\right|=\\
&\mathbb{E}\left|\left(\frac{1}{N}\sum_{i=1}^N\left(S_s^N(i)-\mathbb{E}S_s^N(i)\right)h_1\left(\frac{i}{N}\right)\right)
\left(\frac{1}{N^{0.5}}\sum_{j=1}^N\left(I_s^N(j)-\mathbb{E}I_s^N(j)\right)h_2\left(\frac{j}{N}\right)\right)\right|.
\end{align*}
According to Lemma \ref{lemma 3.1 replacement} and the analysis leading to Lemma \ref{lemma 4.1},
\[
\sup_{0\leq s\leq t}\mathbb{E}\left(\left(\frac{1}{N}\sum_{i=1}^N\left(S_s^N(i)-\mathbb{E}S_s^N(i)\right)h_1\left(\frac{i}{N}\right)\right)^2\right)=O(N^{-1})
\]
and
\[
\sup_{0\leq s\leq t}\mathbb{E}\left(\left(\frac{1}{N^{0.5}}\sum_{j=1}^N\left(I_s^N(j)-\mathbb{E}I_s^N(j)\right)h_2\left(\frac{j}{N}\right)\right)^2\right)=O(1).
\]
By Cauchy-Schwarz inequality, for $h\in \mathcal{H}$,
\begin{equation}\label{equ 5.14}
\sup_{0\leq s\leq t}\mathbb{E}\left|\frac{1}{N^{1.5}}\sum_{i=1}^N\sum_{j=1}^N
\left(S_s^N(i)-\mathbb{E}S_s^N(i)\right)\left(I_s^N(j)-\mathbb{E}I_s^N(j)\right)h\left(\frac{i}{N}, \frac{j}{N}\right)\right|=O(N^{-1/2}).
\end{equation}
As a result, for $h\in \mathcal{H}$, Equation \eqref{equ 5.13} follows from Markov's inequality.

For general $h$, since ${\rm span}(\mathcal{H})$ is dense in $C\left([0, 1]\times[0, 1]\right)$, for any $\varepsilon>0$ there exists $h_\epsilon\in {\rm span}(\mathcal{H})$ such that $\|h_\epsilon-h\|_\infty\leq \epsilon$. According to Lemma \ref{lemma 3.1 replacement}, there exists $C_8<+\infty$ independent of $N, h, h_\epsilon$ such that
\begin{align*}
&\mathbb{E}\left(\left(\frac{1}{N^{1.5}}\sum_{i=1}^N\sum_{j=1}^N\left(S_s^N(i)-\mathbb{E}S_s^N(i)\right)\left(I_s^N(j)-\mathbb{E}I_s^N(j)\right)(h-h_\epsilon)\left(\frac{i}{N}, \frac{j}{N}\right)\right)^2\right)\\
&\leq C_8\|h-h_\epsilon\|^2\leq C_8\varepsilon^2
\end{align*}
for all $0\leq s\leq T$. Therefore, by Equation \eqref{equ 5.14} and H\"{o}lder's inequality,
\begin{align*}
&\sup_{0\leq s\leq t}\mathbb{E}\left|\frac{1}{N^{1.5}}\sum_{i=1}^N\sum_{j=1}^N
\left(S_s^N(i)-\mathbb{E}S_s^N(i)\right)\left(I_s^N(j)-\mathbb{E}I_s^N(j)\right)h\left(\frac{i}{N}, \frac{j}{N}\right)\right|\\
&\leq \sqrt{C_8}\varepsilon+O(N^{-\frac{1}{2}}).
\end{align*}
Let $N\rightarrow+\infty$, then
\begin{align*}
&\lim_{N\rightarrow+\infty}\sup_{0\leq s\leq t}\mathbb{E}\left|\frac{1}{N^{1.5}}\sum_{i=1}^N\sum_{j=1}^N
\left(S_s^N(i)-\mathbb{E}S_s^N(i)\right)\left(I_s^N(j)-\mathbb{E}I_s^N(j)\right)h\left(\frac{i}{N}, \frac{j}{N}\right)\right|\\
&\leq \sqrt{C_8}\varepsilon
\end{align*}
and hence
\[
\lim_{N\rightarrow+\infty}\sup_{0\leq s\leq t}\mathbb{E}\left|\frac{1}{N^{1.5}}\sum_{i=1}^N\sum_{j=1}^N
\left(S_s^N(i)-\mathbb{E}S_s^N(i)\right)\left(I_s^N(j)-\mathbb{E}I_s^N(j)\right)h\left(\frac{i}{N}, \frac{j}{N}\right)\right|=0
\]
since $\varepsilon$ is arbitrary. Consequently, for general $h\in \mathcal{H}$, Equation \eqref{equ 5.13} follows from Markov's inequality.

\qed

At last, we prove Theorem \ref{theorem 2.2 fluctuation}.

\proof[Proof of Theorem \ref{theorem 2.2 fluctuation}]

By Lemma \ref{lemma 5.1}, any subsequence
\[
\left\{\left(\eta_t^{N_j}, \beta_t^{N_j}\right):~0\leq t\leq T\right\}_{j\geq 1}
\]
of
$
\left\{\left(\eta_t^N, \beta_t^N\right):~0\leq t\leq T\right\}_{N\geq 1}
$
has a weakly convergent subsequence
\[\left\{\left(\eta_t^{N_{j_k}}, \beta_t^{N_{j_k}}\right):~0\leq t\leq T\right\}_{k\geq 1}.\] Let $\left(\hat{\eta}_t, \hat{\beta}_t\right)_{0\leq t\leq T}$ be the weak limit of
$\left\{\left(\eta_t^{N_{j_k}}, \beta_t^{N_{j_k}}\right):~0\leq t\leq T\right\}_{k\geq 1}$,
then we only need to show that $\left(\hat{\eta}_t, \hat{\beta}_t\right)_{0\leq t\leq T}$ is the solution to Equation \eqref{equ 2.2 generalized O-U process} with the required initial distribution given in Theorem \ref{theorem 2.2 fluctuation}. For simplicity, from now on we write $\left\{\left(\eta_t^{N_{j_k}}, \beta_t^{N_{j_k}}\right):~0\leq t\leq T\right\}_{k\geq 1}$ as $\left\{\left(\eta_t^N, \beta_t^N\right):~0\leq t\leq T\right\}_{N\geq 1} $.

The fact that $\hat{\beta}_0=-\hat{\eta}_0$ and $\hat{\eta}_0(f)$ follows $\mathbb{N}(0, \int_0^1f^2(u)\phi(u)(1-\phi(u))du)$ for all $f\in C[0, 1]$ follows from Assumption (A). The fact that $\left(\hat{\eta}_t(f), \hat{\beta}_t(f)\right)_{0\leq t\leq T}$ is continuous for any $f\in C[0, 1]$ follows from the fact that
\[
|\eta_t^N(f)-\eta_{t-}^N(f)|, |\beta_t^N(f)-\beta^N_{t-}(f)|\leq \frac{\|f\|_\infty}{\sqrt{N}},
\]
where $t-$ is the moment just before $t$. Therefore, we only need to show that $\{\mathcal{M}_{t,\hat{\eta}, \hat{\beta}}(G, f, g)\}_{0\leq t\leq T}$ is a martingale for any $G\in C_c^\infty(\mathbb{R}^2)$ and $f,g\in C[0, 1]$, where $\mathcal{M}_{t,\hat{\eta}, \hat{\beta}}(G, f, g)$ is defined as
\begin{align*}
&G(\hat{\eta}_t(f), \hat{\beta}_t(g))-G(\hat{\eta}_0(f), \hat{\beta}_0(g))\\
&\text{~}-\int_0^t\frac{\partial}{\partial u}G(\hat{\eta}_s(f), \hat{\beta}_s(g))\left(\hat{\beta}_s(\mathcal{A}_{1,s}f)+\hat{\eta}_s(\mathcal{A}_{0,s}f-\mathcal{A}_\psi f)\right)ds\\
&\text{~}+\int_0^t\frac{\partial}{\partial v}G(\hat{\eta}_s(f), \hat{\beta}_s(g))\left(\hat{\beta}_s(\mathcal{A}_{1,s}g)+\hat{\eta}_s(\mathcal{A}_{0, s}g)\right)ds\\
&\text{~}-\frac{1}{2}\int_0^t\frac{\partial^2}{\partial u^2}G(\hat{\eta}_s(f), \hat{\beta}_s(g))\\
&\text{\quad\quad\quad}\times\left(\int_0^1\rho_1(s,u)\psi(u)f^2(u)du+\int_0^1\int_0^1\lambda(u,v)\rho_0(s,u)\rho_1(s,v)f^2(u)dudv\right)ds\\
&\text{~}-\frac{1}{2}\int_0^t\int_0^1\int_0^1\frac{\partial^2}{\partial v^2}G(\hat{\eta}_s(f), \hat{\beta}_s(g))\lambda(u,v)\rho_0(s,u)\rho_1(s,v)g^2(u)dudvds\\
&\text{~}+\int_0^t\int_0^1\int_0^1\frac{\partial^2}{\partial uv}G(\hat{\eta}_s(f), \hat{\beta}_s(g))\lambda(u,v)\rho_0(s,u)\rho_1(s,v)f(u)g(u)dudvds.
\end{align*}
For all $N\geq 1$, let
\begin{align*}
&\mathcal{M}^N_{t,\eta, \beta}(G, f, g) \\
&=G(\eta_t^N(f), \beta_t^N(g))-G(\eta_0^N(f), \beta_0^N(g))-\int_0^t(\partial_s+\mathcal{L}^N)G(\eta_s^N(f), \beta_s^N(g))ds,
\end{align*}
then $\{\mathcal{M}_{t,\eta, \beta}^N(G, f, g)\}_{0\leq t\leq T}$ is a martingale according to Dynkin's martingale formula. According to the definition of $\mathcal{L}^N$ and Chapman-Kolmogorov equation, we have
\begin{align*}
&\mathcal{L}^NG\left(\eta_s^N(f), \beta_s^N(g)\right)\\
&=\sum_{i=1}^N\psi\left(\frac{i}{N}\right)I_s^N(i)\left(G\left(\eta_s^N(f)-\frac{f\left(\frac{i}{N}\right)}{\sqrt{N}}, \beta_s^N(g)\right)-G(\eta_s^N(f), \beta_s^N(g))\right)\\
&\text{\quad}+\frac{1}{N}\sum_{i=1}^N\sum_{j\neq i}\lambda\left(\frac{i}{N}, \frac{j}{N}\right)S_s^N(i)I_s^N(j)\\
&\text{\quad\quad\quad\quad}\times\left(G\left(\eta_s^N(f)+\frac{f\left(\frac{i}{N}\right)}{\sqrt{N}}, \beta_s^N(g)-\frac{g\left(\frac{i}{N}\right)}{\sqrt{N}}\right)-G(\eta_s^N(f), \beta_s^N(g))\right)
\end{align*}
and
\begin{align*}
&\partial_sG\left(\eta_s^N(f), \beta_s^N(g)\right)\\
&=\frac{1}{\sqrt{N}}\partial_uG(\eta_s^N(f), \beta_s^N(g))\sum_{i=1}^N\psi\left(\frac{i}{N}\right)\mathbb{E}I_s^N(i)f\left(\frac{i}{N}\right)\\
&\text{\quad}-\frac{1}{N^{1.5}}\partial_uG(\eta_s^N(f), \beta_s^N(g))\sum_{i=1}^N\sum_{j\neq i}\lambda\left(\frac{i}{N}, \frac{j}{N}\right)\mathbb{E}\left(S_s^N(i)I_s^N(j)\right)f\left(\frac{i}{N}\right)\\
&\text{\quad}+\frac{1}{N^{1.5}}\partial_vG(\eta_s^N(f), \beta_s^N(g))\sum_{i=1}^N\sum_{j\neq i}\lambda\left(\frac{i}{N}, \frac{j}{N}\right)\mathbb{E}\left(S_s^N(i)I_s^N(j)\right)g\left(\frac{i}{N}\right).
\end{align*}
Then, by Taylor's expansion up to the second order,
\begin{align*}
&(\partial_s+\mathcal{L}^N)G(\eta_s^N(f), \beta_s^N(g)) \\
&=\partial_uG(\eta_s^N(f), \beta_s^N(g))\Bigg(-\eta_s^N(\mathcal{A}_\psi f)\\
&\text{\quad\quad\quad}+\frac{1}{N^{1.5}}\sum_{i=1}^N\sum_{j\neq i}\lambda\left(\frac{i}{N}, \frac{j}{N}\right)\left(S_s^N(i)I_s^N(j)-\mathbb{E}\left(S_s^N(i)I_s^N(j)\right)\right)f\left(\frac{i}{N}\right)\Bigg)\\
&\text{\quad}-\partial_vG(\eta_s^N(f), \beta_s^N(g))\times\\
&\text{\quad\quad\quad}\left(\frac{1}{N^{1.5}}\sum_{i=1}^N\sum_{j\neq i}\lambda\left(\frac{i}{N}, \frac{j}{N}\right)\left(S_s^N(i)I_s^N(j)-\mathbb{E}\left(S_s^N(i)I_s^N(j)\right)\right)g\left(\frac{i}{N}\right)\right)\\
&\text{\quad}+\frac{1}{2}\partial_{uu}^2G(\eta_s^N(f), \beta_s^N(g))\Bigg(\frac{1}{N}\sum_{i=1}^N\psi\left(\frac{i}{N}\right)I_s^N(i)f^2\left(\frac{i}{N}\right)\\
&\text{\quad\quad\quad}+\frac{1}{N^2}\sum_{i=1}^N\sum_{j\neq i}\lambda\left(\frac{i}{N}, \frac{j}{N}\right)S_s^N(i)I_s^N(j)f^2\left(\frac{i}{N}\right)\Bigg)\\
&\text{\quad}-\frac{\partial_{uv}^2G(\eta_s^N(f), \beta_s^N(g))}{N^2}\sum_{i=1}^N\sum_{j\neq i}\lambda\left(\frac{i}{N}, \frac{j}{N}\right)S_s^N(i)I_s^N(j)f\left(\frac{i}{N}\right)g\left(\frac{i}{N}\right)\\
&\text{\quad}+\frac{1}{2}\frac{\partial_{vv}^2G(\eta_s^N(f), \beta_s^N(g))}{N^2}\sum_{i=1}^N\sum_{j\neq i}\lambda\left(\frac{i}{N}, \frac{j}{N}\right)S_s^N(i)I_s^N(j)g^2\left(\frac{i}{N}\right)+\varepsilon_{s,1}^N,
\end{align*}
where $\varepsilon_{s,1}^N$ is the Lagrange remainder such that
\[
\|\varepsilon_{s, 1}^N\|\leq \frac{3\|G^{\prime\prime\prime}\|_\infty\left(\|f\|_\infty^3+\|g\|_\infty^3\right)\left(\|\psi\|_\infty+\|\lambda\|_\infty\right)}{2\sqrt{N}}.
\]
According to the above equation, Lemma \ref{lemma 3.1 replacement} and Cauchy-Schwarz inequality, there exists $C_6<\infty$ independent of $N$ such that
\[
\sup_{N\geq 4, 0\leq t\leq T}\mathbb{E}\left(\left(\mathcal{M}^N_{t,\eta, \beta}(G, f, g)\right)^2\right)\leq C_6.
\]
Hence, for any $0\leq t\leq T$, $\{\mathcal{M}^N_{t, \eta, \beta}(G, f, g)\}_{N\geq 4}$ are uniformly integrable. Therefore, according to Theorem 5.3 of \cite{Whitt2007}, we only need to show that $\mathcal{M}^N_{t, \eta, \beta}(G, f, g)$ converges weakly to $\mathcal{M}_{t,\hat{\eta}, \hat{\beta}}(G, f, g)$ as $N\rightarrow+\infty$ to complete this proof.

According to Lemma \ref{lemma 3.1 replacement} and Cauchy-Schwarz inequality, random variables $I_s^N(i)$ and $S_s^N(i)I_s^N(j)$ in the last three terms of the expression of $(\partial_s+\mathcal{L}^N)G(\eta_s^N(f), \beta_s^N(g))$ can be replaced by their means with small errors. Similarly, according to Lemmas \ref{lemma 3.1 replacement} and \ref{lemma 5.2}, the random variable
\[
S_s^N(i)I_s^N(j)-\mathbb{E}\left(S_s^N(i)I_s^N(j)\right)
\]
in the first two terms of the expression of $(\partial_s+\mathcal{L}^N)G(\eta_s^N(f), \beta_s^N(g)) $ can be replaced by
\[
(\mathbb{E}S_s^N(i))\left(I_s^N(j)-\mathbb{E}I_s^N(j)\right)+(\mathbb{E}I_s^N(j))\left(S_s^N(i)-\mathbb{E}S_s^N(i)\right)
\]
with a small error since
\begin{align*}
&S_s^N(i)I_s^N(j)-\mathbb{E}\left(S_s^N(i)I_s^N(j)\right)\\
&\text{\quad}-\left((\mathbb{E}S_s^N(i))\left(I_s^N(j)-\mathbb{E}I_s^N(j)\right)+(\mathbb{E}I_s^N(j))\left(S_s^N(i)-\mathbb{E}S_s^N(i)\right)\right)\\
&=\left(S_s^N(i)-\mathbb{E}S_s^N(i)\right)\left(I_s^N(j)-\mathbb{E}I_s^N(j)\right)-\left(\mathbb{E}\left(S_s^N(i)I_s^N(j)\right)
-\mathbb{E}S_s^N(i)\mathbb{E}I_s^N(j)\right).
\end{align*}

In detail,
\begin{align*}
&(\partial_s+\mathcal{L}^N)G(\eta_s^N(f), \beta_s^N(g)) \\
&=\varepsilon_{s,1}^N+\varepsilon_{s,2}^N+\partial_uG(\eta_s^N(f), \beta_s^N(g))\Big(-\eta_s^N(\mathcal{A}_\psi f)
+\eta_s^N(\mathcal{A}_{0,s}^Nf)+\beta_s^N(\mathcal{A}_{1,s}^Nf)\Big)\\
&\text{\quad}-\partial_vG(\eta_s^N(f), \beta_s^N(g))
\left(\beta_s^N(\mathcal{A}_{1,s}^Ng)+\eta_s^N(\mathcal{A}_{0,s}^Ng)\right)\\
&\text{\quad}+\frac{1}{2}\partial_{uu}^2G(\eta_s^N(f), \beta_s^N(g))\Bigg(\frac{1}{N}\sum_{i=1}^N\psi\left(\frac{i}{N}\right)\mathbb{E}I_s^N(i)f^2\left(\frac{i}{N}\right)\\
&\text{\quad\quad\quad}+\frac{1}{N^2}\sum_{i=1}^N\sum_{j\neq i}\lambda\left(\frac{i}{N}, \frac{j}{N}\right)\mathbb{E}\left(S_s^N(i)I_s^N(j)\right)f^2\left(\frac{i}{N}\right)\Bigg)\\
&\text{\quad}-\frac{\partial_{uv}^2G(\eta_s^N(f), \beta_s^N(g))}{N^2}\sum_{i=1}^N\sum_{j\neq i}\lambda\left(\frac{i}{N}, \frac{j}{N}\right)\mathbb{E}\left(S_s^N(i)I_s^N(j)\right)f\left(\frac{i}{N}\right)g\left(\frac{i}{N}\right)\\
&\text{\quad}+\frac{1}{2}\frac{\partial_{vv}^2G(\eta_s^N(f), \beta_s^N(g))}{N^2}\sum_{i=1}^N\sum_{j\neq i}\lambda\left(\frac{i}{N}, \frac{j}{N}\right)\mathbb{E}\left(S_s^N(i)I_s^N(j)\right)g^2\left(\frac{i}{N}\right),
\end{align*}
where
\[
\lim_{N\rightarrow+\infty}\int_0^t\varepsilon_{s,1}^Nds=0
\]
in probability,
\[
\mathbb{E}\left(\left(\int_0^t \varepsilon_{2,s}^Nds\right)^2\right)\leq \frac{C_7}{N}
\]
for some $C_7<+\infty$ independent of $N$,
\[
\mathcal{A}_{0,s}^Nf(u)=\frac{1}{N}\sum_{j=1}^N\lambda\left(\frac{j}{N}, u\right)\mathbb{E}S_s^N(j)f\left(\frac{j}{N}\right)
\]
and
\[
\mathcal{A}_{1,s}^Nf(u)=f(u)\frac{1}{N}\sum_{j=1}^N\lambda\left(u, \frac{j}{N}\right)\mathbb{E}I_s^N(j).
\]
According to Equation \eqref{equ 4.1}, Lemma \ref{lemma 3.1 replacement} and Chebyshev's inequality,
\[
\lim_{N\rightarrow+\infty}\int_0^t\eta_s^N\left(\mathcal{A}_{0,s}(f)-\mathcal{A}_{0,s}^N(f)\right)ds=0
\]
and
\[
\lim_{N\rightarrow+\infty}\int_0^t\beta_s^N\left(\mathcal{A}_{1,s}(f)-\mathcal{A}_{1,s}^N(f)\right)ds=0
\]
in probability. By Lemma \ref{lemma 3.1 replacement}, $\mathbb{E}\left(S_s^N(i)I_s^N(j)\right)$ in the last three terms of the above expression of $(\partial_s+\mathcal{L}^N)G(\eta_s^N(f), \beta_s^N(g))$ can be replaced by $\mathbb{E}S_s^N(i)\mathbb{E}I_s^N(j)$ with small error. Furthermore, by Equation \eqref{equ 4.1}, in the new expression of $(\partial_s+\mathcal{L}^N)G(\eta_s^N(f), \beta_s^N(g))$ we can replace $\frac{1}{N}\sum_{i=1}^N\psi\left(\frac{i}{N}\right)\mathbb{E}I_s^N(i)f^2\left(\frac{i}{N}\right)$ by
\[
\int_0^1\psi(u)\rho_1(s,u)f^2(u)du
\]
and replace any term with form
$\frac{1}{N^2}\sum_{i=1}^N\sum_{j\neq i}\lambda\left(\frac{i}{N}, \frac{j}{N}\right)\mathbb{E}S_s^N(i)\mathbb{E}I_s^N(j)h\left(\frac{i}{N}, \frac{j}{N}\right)$
by
\[
\int_0^1\int_0^1\lambda(u,v)\rho_0(s,u)\rho_1(s,v)h(u,v)dudv
\]
with small errors. In conclusion,
\begin{align*}
&(\partial_s+\mathcal{L}^N)G(\eta_s^N(f), \beta_s^N(g)) \\
&=\varepsilon_{s,4}^N+\partial_uG(\eta_s^N(f), \beta_s^N(g))\Big(-\eta_s^N(\mathcal{A}_\psi f)
+\eta_s^N(\mathcal{A}_{0,s}f)+\beta_s^N(\mathcal{A}_{1,s}f)\Big)\\
&\text{\quad}-\partial_vG(\eta_s^N(f), \beta_s^N(g))
\left(\beta_s^N(\mathcal{A}_{1,s}g)+\eta_s^N(\mathcal{A}_{0,s}g)\right)\\
&\text{\quad}+\frac{1}{2}\partial_{uu}^2G(\eta_s^N(f), \beta_s^N(g))\Bigg(\int_0^1\psi(u)\rho_1(s, u)f^2(u)du\\
&\text{\quad\quad\quad}+\int_0^1\int_0^1\lambda(u,v)\rho_0(s,u)\rho_1(s,v)f^2(u)dudv\Bigg)\\
&\text{\quad}-\partial_{uv}^2G(\eta_s^N(f), \beta_s^N(g))\int_0^1\int_0^1\lambda(u, v)\rho_0(s,u)\rho_1(s,v)f(u)g(u)dudv\\
&\text{\quad}+\frac{1}{2}\partial_{vv}^2G(\eta_s^N(f), \beta_s^N(g))\int_0^1\int_0^1\lambda(u,v)\rho_0(s,u)\rho_1(s,v)g^2(u)dudv,
\end{align*}
where
\[
\lim_{N\rightarrow+\infty}\int_0^t \varepsilon_{s, 4}^Nds=0
\]
in probability. Since $\left\{\left(\hat{\eta}_t, \hat{\beta}_t\right):~0\leq t\leq T\right\}$ is a weak limit of
\[
\left\{\left(\eta_t^N, \beta_t^N\right):~0\leq t\leq T\right\}_{N\geq 1},
\]
let $N\rightarrow+\infty$ in the above expression of $(\partial_s+\mathcal{L}^N)G(\eta_s^N(f), \beta_s^N(g)) $, then we have $\mathcal{M}^N_{t, \eta, \beta}(G, f, g)$ converges weakly to $\mathcal{M}_{t,\hat{\eta}, \hat{\beta}}(G, f, g)$ as $N\rightarrow+\infty$ and the proof is complete.

\qed

\quad

\textbf{Acknowledgments.} The authors are grateful to the financial
support from the National Natural Science Foundation of China with
grant number 11501542.

{}
\end{document}